\pdfoutput=1
\documentclass[11pt,francais]{smfart}
\usepackage{Publier}

\begin{document}

\title[Publier sous l'Occupation]{Publier sous l'Occupation~I. Autour~du~cas~de~Jacques~Feldbau et~de~l'Acad\'emie~des~sciences}
\date{\today}
\author{Mich\`ele Audin}
\thanks{Mich\`ele Audin, Institut de Recherche math\'ematique avanc\'ee, Universit\'e Louis Pasteur et CNRS, 7 rue Ren\'e Descartes, 67084 Strasbourg Cedex, France}
\address{Institut de Recherche Math\'ematique Avanc\'ee\\
Universit\'e Louis Pasteur et CNRS\\
7 rue Ren\'e Descartes\\
67084 Strasbourg cedex\\
France}
\email{Michele.Audin@math.u-strasbg.fr}
\urladdr{http://www-irma.u-strasbg.fr/~maudin}


\maketitle

\begin{abstract}
C'est un article sur les publications math\'ematiques pendant l'Occupation (1940--44). \`A travers les cas de quatre d'entre eux, et surtout de celui de Jacques Feldbau (un des fondateurs de la th\'eorie des fibr\'es, mort en d\'eportation), nous \'etudions la fa\c con dont la censure a frapp\'e les math\'ematiciens fran\c cais d\'efinis comme juifs par le \og Statut des juifs\fg d'octobre 1940 et les strat\'egies de publication que ceux-ci ont alors utilis\'ees (pseudonymes, plis cachet\'es, journaux provinciaux...) La mani\`ere dont les \og lois en vigueur\fg ont \'et\'e (ou n'ont pas \'et\'e) discut\'ees et appliqu\'ees \`a l'Acad\'emie des sciences est \'egalement \'etudi\'ee.
\end{abstract}

\begin{altabstract}
This is an article on mathematical publishing during the German Occupation of France (1940--44). Looking at the cases of four of them and especially at the case of Jacques Feldbau (one of the founders of the theory of fibre bundles, dead in deportation), we investigate the way censorship struck the French mathematicians who were declared jewish by the \og Statut des juifs\fg of october 1940, and the strategies these mathematicians then developed (fake names, selled envelopes, provincial journals...). The way the Vichy laws have been (or have not been) discussed and applied at the Acad\'emie des sciences is investigated as well.
\end{altabstract}

\section*{Introduction}
L'objet de cet article est d'\'etudier, autour du cas du topologue Jacques Feldbau, la mani\`ere dont les math\'ematiciens fran\c cais juifs (c'est-\`a-dire d\'efinis comme tels par le \og Statut des juifs\fg du 3 octobre 1940) ont pu --- ou n'ont pas pu --- publier les r\'esultats de leurs recherches dans les journaux sp\'ecialis\'es pendant la p\'eriode de l'Occupation allemande (1940--1944).

\subsection*{Les contextes}

Le contexte pour les math\'ematiques est celui de l'\'elaboration dans les ann\'ees 1930 et 1940, des fondements d'un groupe de sous-disciplines que l'on appelle aujourd'hui la topologie (g\'en\'erale, alg\'ebrique, diff\'erentielle) et la g\'eom\'etrie (diff\'erentielle) et en particulier de la participation \`a cette \'elaboration de math\'{e}maticiens fran\c cais (de la descendance d'\'Elie Cartan), Ehresmann et Feldbau notamment.

\medskip
Le contexte historique g\'en\'eral, dont ce contexte math\'ematique est indissociable, est celui de l'Occupation allemande et du r\'egime de Vichy. Rappelons que le territoire fran\c cais est en majorit\'e occup\'e par les troupes allemandes \`a partir de l'armistice du 22 juin 1940 (il le sera compl\`etement apr\`es le 11 novembre 1942); le pays (priv\'e de l'Alsace et de la Moselle) est administr\'e depuis Vichy par le gouvernement de l'\og\'Etat fran\c cais\fg dirig\'e par Philippe P\'etain --- l'Allemagne exer\c cant de surcro\^it ses droits de puissance occupante dans la zone occup\'ee.

On le sait, cet \og\'Etat fran\c cais\fg a tr\`es rapidement et largement anticip\'{e} et devanc\'e les d\'esirs des occupants, notamment, pour ce qui nous concerne ici, en promulguant, d\`es le 3 octobre 1940, une s\'erie de d\'ecrets rassembl\'es sous le titre de \og Statut des juifs\fg. Le cadre g\'en\'eral des effets de la politique de Vichy,  et en particulier de ces d\'ecrets, sur l'Universit\'e est bien d\'ecrit et \'etudi\'e dans les travaux pr\'ecurseurs de Claude Singer, notamment dans son livre~\citepyear{singer92} et dans son article~\citepyear{singervichy}.

On sait aussi que cet antis\'emitisme fran\c cais officiel s'est d\'evelopp\'e et renforc\'e, des d\'ecrets du 3 octobre 1940 \`a celui du 6 juin 1942, en passant par la loi du 2 juin 1941 et l'\'etablissement d'un fichier des juifs, allant d'une logique d'exclusion des juifs \`a la logique d'extermination\footnote{Comme le dit l'historien Denis Peschanski dans l'introduction de~\cite{drancy}.} qui lui a fait rafler les Fran\c cais juifs apr\`es les juifs d'origines \'etrang\`ere pour les regrouper dans des camps comme celui de Drancy avant de les acheminer vers les camps d'extermination nazis --- ce sera notamment le sort de Jacques Feldbau, math\'ematicien fran\c cais juif mort en d\'eportation.

En ce qui concerne les publications scientifiques, le statut des juifs du 3 octobre 1940, dans son article 5, stipule:
\begin{quote}
Les juifs ne pourront, sans condition ni r\'eserve, exercer l'une quelconque des professions suivantes: 
Directeurs, g\'erants, r\'edacteurs de journaux, revues, agences ou p\'eriodiques, \`a l'exception de publications de caract\`ere strictement scientifique. [...]
\end{quote}
La m\^eme exception (qui semble tol\'erer que les scientifiques juifs publient leurs articles) figure encore dans le statut modifi\'e du 2 juin 1941 (toujours l'article~5)\footnote{\emph{Journal officiel}, 14 juin 1941, p.~2475.}:
\begin{quote}
Sont interdites aux juifs les professions ci-apr\`es: 
Banquier, changeur, d\'emarcheur; [...]

\'Editeur, directeur, g\'erant, administrateur, r\'edacteur, m\^eme au titre de correspondant local, de journaux ou d'\'ecrits p\'eriodiques, \`a l'exception des publications de caract\`ere strictement scientifique ou confessionnel; [...]
\end{quote}

Il semblerait donc qu'un scientifique, m\^eme r\'eput\'e juif, puisse continuer, d'apr\`es la loi fran\c caise, \`a publier. La r\'ealit\'e, nous allons le voir, est assez diff\'erente. Beaucoup de journaux scientifiques sont publi\'es \`a Paris, c'est aussi l\`a que si\`ege l'Acad\'emie des sciences, en zone occup\'ee, donc. Et il y a une censure allemande pendant l'Occupation, il y a aussi des interdits sur les publications et, en coh\'erence avec la politique de l'\og \'Etat fran\c cais\fg, devan\c cant les demandes allemandes comme pour le statut des juifs, il y a des scientifiques fran\c cais qui participent \`a cette censure. 

Une note apparaissant dans la th\`ese~\citepyear{pinaultthese} de Michel Pinault (mais pas dans le livre~\citepyear{pinaultlivre} qui en est tir\'e) \'evoque le r\^ole jou\'e par Ernest Fourneau\footnote{Le chimiste Ernest Fourneau (1872--1949) a \'et\'e pr\'esident du Comit\'e consultatif de la litt\'erature et de la documentation scientifiques du Groupement de la presse p\'eriodique de 1942 \`a 1944. Il a eu quelques ennuis (de courte dur\'ee) \`a la Lib\'eration.}, un membre de l'Acad\'emie de m\'edecine, qui dirigeait l'Institut Pasteur et \`a qui les autorit\'es allemandes avaient donn\'e tout pouvoir sur les publications scientifiques.
\begin{quote}
Les conditions dans lesquelles les Allemands confient \`a Fourneau la responsabilit\'e de donner, en leur nom, l'autorisation de para\^itre aux publications scientifiques, ne sont pas claires. Elles r\'esultent probablement de Fourneau lui-m\^eme qui leur aurait offert ses services.~\cite{pinaultthese}
\end{quote}
Et il cite une lettre du \textsc{mbf}\footnote{Le sigle \textsc{mbf} d\'esigne le \emph{Milit\"arbefehlshaber in Frankreich}, l'autorit\'e d'occupation militaire.} adress\'ee \`a Fourneau, le 28 novembre 1940:
\begin{quote}
\`A la suite de notre entretien du 26 courant, j'ai l'honneur de vous informer de ce qui suit: 1) Pour le r\`{e}glement int\'erieur de la Vie universitaire fran\c caise, seules les lois de votre Gouvernement sont valables. Il en est donc ainsi \'egalement pour la question juive. 2) Il est \'evident que, dans le sens de la collaboration sinc\`erement recherch\'ee par nous, nous serions heureux si les savants fran\c cais allaient d\'ej\`a plus loin que les premi\`eres lois de Vichy sur les juifs et s'attaquaient au probl\`eme d'\'epurer leur profession des juifs qui en font partie. 3) Sans vouloir vous donner des directives \`a ce sujet, nous consid\'ererions avec plaisir, et comme allant de soi, qu'aucun juif ne nous soit pr\'esent\'e dans les organes repr\'esentatifs de votre profession qui, pour une raison quelconque, doivent collaborer avec les autorit\'es d'occupation, c'est-\`a-dire dans les comit\'es qui, par exemple, traitent avec nous des questions relatives \`a la censure des publications et brochures scientifiques. De m\^eme nous pensons qu'il va de soi que, dans vos publications et brochures scientifiques, aucun juif ne paraisse sur la couverture. Vous jugerez naturellement vous-m\^eme pendant combien de temps vous croyez ne pas pouvoir vous passer de la collaboration de coll\`egues juifs, mais il nous semble qu'une adaptation la plus rapide possible \`a l'exemple allemand est parfaitement possible d'apr\`es les exp\'eriences faites chez nous. [Lettre du \textsc{mbf}--\emph{Propaganda Abteilung}, R\'ef. \emph{Schrifften} Pgb Nr 01166/40, adress\'ee \`a \og monsieur le professeur Fourneau\fg et sign\'ee: Schulz, le 28 novembre 1940), cit\'ee dans~\cite{pinaultthese}.]
\end{quote}

\subsection*{La probl\'ematique}

La probl\'ematique envisag\'ee dans cet article est double:
\begin{itemize}
\item d'une part, quel a \'et\'e l'effet de ce contexte politique sur l'histoire de la mise en place des outils fondamentaux de la topologie,
\item d'autre part, comment les institutions se sont accommod\'{e}es\footnote{Voir la note~\ref{accommodement}.} des lois en vigueur pour emp\^echer les math\'ematiciens relevant des diff\'erents statuts des juifs de publier, et quelles strat\'egies ceux-ci ont d\'evelopp\'ees pour faire conna\^itre leurs travaux.
\end{itemize}

\subsection*{Les math\'ematiciens concern\'es}

Plusieurs math\'{e}maticiens menac\'es par les dispositifs antis\'emites de Vichy que nous venons d'\'evoquer ont quitt\'e la France et pass\'e la plus grande partie de la p\'eriode de l'Occupation aux \'Etats-Unis (Hadamard, Mandelbrojt, Weil...) et ils ont naturellement publi\'e dans les journaux sp\'ecialis\'es am\'ericains. Il n'en sera donc pas question ici.

Mon point de d\'epart a \'et\'e le travail math\'ematique du topologue Jacques Feldbau. Ce travail, dont le contenu est important pour l'histoire de la topologie, a \'et\'e publi\'e d'une mani\`ere totalement anormale, une anormalit\'e li\'ee \`a plusieurs titres aux questions pos\'ees ci-dessus: une notion fondamentale de la th\'eorie (celle de fibr\'e associ\'e) est publi\'ee sans qu'il apparaisse en \^etre auteur, une autre notion importante (le produit de Whitehead) para\^it dans un article treize ans apr\`es sa mort --- et plus de quinze ans apr\`es avoir \'et\'e d\'{e}fini par Whitehead ---, deux de ses articles sont publi\'es sous un pseudonyme...

Il \'etait tentant d'essayer de savoir si d'autres math\'ematiciens avaient v\'ecu des exp\'eriences analogues. J'ai recens\'e les cas de Bloch, Schwartz et Paul L\'evy.

Andr\'e Bloch a lui aussi publi\'e sous (deux) pseudonyme(s). La diserte autobiographie~\cite{schwartz}\footnote{Puisqu'il est question d'autobiographies, t\'emoignages de leurs auteurs sur cette \'epoque, on sait que les souvenirs~\cite{weilsouvenirs} de Weil ne donnent aucune information pertinente pour le sujet \'etudi\'e ici. Le dernier article publi\'e par Weil en France est une c\'el\`ebre note aux \emph{Comptes rendus}~\cite{weil40} (voir~\cite[p.~546]{weilOC1} et~\cite[p.~153]{weilsouvenirs}) qui date du 22 avril 1940, avant l'Occupation, donc. Apr\`es diverses tribulations (dont il rend compte dans son \og ballet-bouffe\fg), il quitte la France pour les \'Etats-Unis d\`es janvier 1941~\cite[p.~180]{weilsouvenirs}.} de Laurent Schwartz est malheureusement tr\`es impr\'ecise sur la mani\`ere dont son auteur a publi\'e pendant l'Occupation, nous verrons qu'il a trouv\'e un \'editeur et m\^eme un journal pour l'accueillir. La correspondance~\cite{LevyFrechet} entre Paul L\'evy et Maurice Fr\'echet donne, elle, de pr\'ecieuses informations sur le contexte tout en montrant \`a quel point Paul L\'evy\footnote{Une personnalit\'e attachante dont la biographie est bien document\'ee (voir notamment~\cite{SchwartzLevy,LevyFrechet,Levy70}).} \'etait attentif aux questions de priorit\'e. On le sait, il publiait vite (souvent trop vite) et beaucoup, sa strat\'egie en mati\`ere de publications \'etait donc \`a \'etudier. J'\'evoquerai aussi bri\`evement le cas de F\'elix Pollaczek, qui \'etait d\'ej\`a install\'e en France et y a v\'ecu cette p\'eriode mais sur lequel j'ai tr\`es peu d'information (voir~\cite{CohenPollaczek}).

\medskip
Feldbau, Bloch, Schwartz, L\'evy et Pollaczek, pourquoi s'arr\^eter l\`a? La r\'eponse est simple: \`a part Marie-H\'el\`ene Schwartz (qui n'a publi\'e qu'une note aux \emph{Comptes rendus} pendant l'Occupation), je connais peu d'autres exemples de math\'ematiciens vivant en France, atteints par les lois antis\'emites, et publiant pendant cette p\'eriode\footnote{Il est impossible de savoir si ceux qui n'ont pas publi\'e pendant cette p\'eriode ne travaillaient plus, pour une raison ou pour une autre, ou ont appliqu\'e une auto-censure. On peut penser par exemple \`a Ayzyk Gorny, math\'ematicien d'origine polonaise et \'el\`eve de Mandelbrojt, qui s'est trouv\'e \`a Clermont-Ferrand comme Schwartz et Feldbau, dans des conditions d'existence pr\'ecaires (d'apr\`es~\cite[p~157]{schwartz}) et qui a \'et\'e d\'eport\'e, comme \og juif \'etranger\fg, dans le convoi 37 le 25~septembre~1942 et est mort \`a Auschwitz (voir~\cite{coutyglaeserperol}), mais ses derniers articles publi\'es l'ont \'et\'e en 1939, assez longtemps donc avant l'Occupation et le Statut des juifs.}. Il n'est d'ailleurs pas tr\`es \'etonnant que cette courte liste contienne deux math\'ematiciens en tout d\'ebut de carri\`ere et donc inconnus (Feldbau et Schwartz) ainsi qu'un math\'ematicien enferm\'e dans un h\^opital psychiatrique (Bloch): il n'y aurait eu aucune raison, respectivement aucun moyen, de les inviter aux \'Etats-Unis. Quant \`a Paul L\'evy, sa correspondance montre \`a quel point il \'etait peu conscient des dangers qu'il courait.

\subsubsection*{Une pr\'ecaution} 
Les dispositifs antis\'emites dont les math\'ema\-ticiens mentionn\'es ici ont \'et\'e victimes \`a des degr\'es divers excluaient des Fran\c cais \emph{qu'eux-m\^emes avaient d\'efinis comme juifs}\footnote{\emph{Journal Officiel}, 18 octobre 1940, p.~5323. L'article premier stipule: \og Est regard\'e comme juif, pour l'application de la pr\'esente loi, toute personne issue de trois grands-parents de race juive ou de deux grands-parents de la m\^eme race, si son conjoint lui-m\^{e}me est juif\fg.\label{notestatut} Une d\'efinition dont les diverses absurdit\'es n'auront \'echapp\'e \`a aucun math\'ematicien.}. C'est le cas des cinq math\'ematiciens dont il est question ici (comme Schwartz l'explique tr\`es bien dans son livre~\citepyear{schwartz}). Pour cette raison ou pour une autre, les institutions ou les coll\`egues r\'edacteurs qui d\'ecidaient de publier ou de ne pas publier leurs travaux  les consid\'eraient comme juifs. Il n'y avait pas lieu d'entrer ici dans d'autres consid\'erations.

\subsection*{M\'ethodologie}

\`A cette \'epoque et m\^eme si elle n'en a d\'ej\`a plus l'exclusivit\'e, l'Acad\'emie des sciences joue encore, \`a travers ses \emph{Comptes rendus}, un r\^ole important de validation des r\'esultats scientifiques. De plus ses archives constituent une source exceptionnellement riche de renseignements. Cette institution (en tant que telle) a r\'eussi \`a s'accommoder\footnote{J'utilise ici l'utile notion d'\og accommodement\fg due \`a Philippe Burrin~\citepyear{burrin}.\label{accommodement}} des lois en vigueur --- avec une prudente discr\'etion. De m\^eme qu'une \'etude syst\'ematique des archives pertinentes peut permettre d'\'ecrire la \og biographie\fg d'un inconnu destin\'e \`a rester un parfait anonyme (voir~\cite{corbin}), l'\'etude syst\'ematique de celles de l'Acad\'emie des sciences permet de faire une histoire de l'application de ces lois, m\^eme si, par exemple, l'interdiction de publier faite aux auteurs juifs n'y est jamais mentionn\'ee en tant que telle\footnote{Absence de mention ne signifie pas que ces mentions auraient \'et\'e d\'etruites. Il est clair que rien n'a disparu des pochettes de s\'eances, qu'aucune page n'a \'et\'e arrach\'ee au registre des comit\'es secrets. La chose n'a pas \'et\'e mentionn\'ee du tout. Il est vrai aussi que le contenu des dossiers biographiques d\'epend \'evidemment de qui a conserv\'e quoi et \`a quelle date --- il est donc plus al\'eatoire.\label{notedestruction}}. J'ai donc lu \emph{toute} la correspondance re\c cue (en particulier les manuscrits des notes aux \emph{Comptes rendus}) et conserv\'ee dans les pochettes hebdomadaires des s\'eances (un travail \`a la fois fastidieux et passionnant), ainsi que beaucoup d'autres documents disponibles\footnote{Dans un souci de rigueur \'el\'ementaire, je cite pr\'ecis\'ement toutes les sources dans le texte.}. Ajoutons que l'Acad\'emie des sciences re\c coit une correspondance, publie des notes, d'un \'eventail plus large de scientifiques que les seuls math\'ematiciens, nous verrons donc d'autant plus de traces du ph\'enom\`ene \'etudi\'e. Car, nous le verrons, les traces sont t\'enues, mais elles existent. 

\medskip
Il est remarquable que les journaux scientifiques fran\c cais (autres que les \emph{Comptes rendus}) cit\'es ici ne semblent pas avoir conserv\'e d'archives. Il est vrai que la p\'eriode \'etait dangereuse et que beaucoup d'\'ecrits ont \'et\'e d\'etruits pour des raisons de prudence \'el\'ementaire... Les renseignements que j'ajoute ici \`a propos de ces journaux viennent donc exclusivement de la lecture de sources publi\'ees (et cit\'ees, elles aussi, dans le texte). Il convient de les consid\'erer surtout comme des compl\'ements d'information.

\subsection*{Sources}
La source principale de ce travail, outre les articles publi\'es de Jacques Feldbau, d'Andr\'e Bloch, de Paul L\'evy, de Laurent Schwartz, de F\'elix Pollaczek ou d'autres, est l'ensemble des pochettes de s\'eances conserv\'ees aux archives de l'Acad\'emie des sciences, accompagn\'e des registres des comit\'es secrets et des dossiers biographiques de certains acad\'emiciens. 

Les correspondances de Hasse et de S\"uss cit\'ees bri\`evement ici sont conserv\'ees respectivement aux universit\'es de G\"{o}ttingen (\emph{Nachla\ss\ H. Hasse}) et de Freiburg (\emph{Nachla\ss\ W. S\"uss}).

\section{Le cas de Jacques Feldbau (1914--1945)}
Jacques Feldbau est n\'e \`a Strasbourg le 22 octobre 1914, il y a fait ses \'etudes, \`a l'exception d'une ann\'ee d'auditorat libre \`a l'\textsc{ens} pour pr\'eparer l'agr\'egation en 1938. Il a \'et\'e le premier \'el\`eve d'Ehresmann\footnote{Charles Ehresmann (1905--1979) a soutenu sa th\`ese \`a Paris en 1934. Il a ensuite \'et\'e nomm\'e \`a Strasbourg. C'est un des inventeurs de la g\'eom\'etrie diff\'erentielle et un des fondateurs du groupe Bourbaki.}.

\subsection{La liste de publications de Jacques Feldbau}
L'\oe uvre math\'ema\-tique publi\'ee de Jacques Feldbau n'est constitu\'ee que d'une trentaine de pages. Voici une liste chronologique des r\'ef\'erences \`a ses articles dans la bibliographie du pr\'esent texte: \cite{Feldbau39}, \cite{ehresmannfeldbau}, \cite{Ehresmann41}, \cite{Laboureur42}, \cite{Laboureur43}, \cite{Feldbau58}. Elle a \'et\'e \'etablie \`a partir des informations contenues dans \emph{Math. Reviews} (sous les noms d'auteurs Feldbau et Laboureur) et de celles dues \`a Michel Zisman~\cite{zismanjames}. On remarquera dans cette liste un certain nombre d'anomalies:
\begin{itemize}
\item un article, \cite{Ehresmann41}, dont Feldbau n'est pas auteur,
\item deux articles, \cite{Laboureur42,Laboureur43}, publi\'es sous un pseudonyme,
\item les \emph{C. R. Acad. Sci. Paris} remplac\'es par le \emph{Bull. Soc. Math. France},
\item un article, \cite{Feldbau58}, publi\'e plus de treize ans apr\`es la mort de son auteur.
\end{itemize}
Un des buts de cette partie est d'expliquer ces anomalies et en particulier pourquoi cette liste est bien la liste des publications de Jacques Feldbau.

\subsection{La topologie dans les ann\'ees 1930--1950 et les notes de Feldbau}
\subsubsection{Sur la classification des espaces fibr\'es}
La premi\`ere note~\citepyear{Feldbau39} de Jacques Feldbau (pr\'esent\'ee par \'Elie Cartan\footnote{\'Elie Cartan (1869--1951), membre de l'Acad\'emie des sciences depuis 1931, avait dirig\'e la th\`ese d'Ehresmann et c'est donc lui qui, naturellement, a transmis cette note de son \'el\`eve Feldbau.} \`a la s\'eance du 15 mai 1939, mais parue dans le fascicule du 22 mai) participe au tout d\'ebut de la th\'eorie des fibr\'es. Il faut comprendre en effet que c'est \`a cette \'epoque que l'on \emph{inventait} les espaces fibr\'es, et que c'est \`a cette invention que Feldbau contribue dans cette note, ce dont t\'emoigne la premi\`ere phrase:
\begin{quote}
Les espaces fibr\'es ont \'et\'e introduits par M.~Seifert dans le cas de vari\'et\'es \`a 3 dimensions et de fibres circulaires. M.~Whitney\footnote{Herbert Seifert (1907--1996) a soutenu sa th\`ese \emph{Topologie dreidimensionaler gefaserter R\"aume} en 1932 \`a Leipzig. Hassler Whitney (1907--1989) a fait une th\`ese sur la th\'eorie des graphes \`a Harvard en 1932, c'est un des fondateurs de la th\'eorie des vari\'et\'es et de la topologie diff\'erentielle.} a \'etudi\'e certains espaces fibr\'es par des sph\`eres. Nous nous proposons d'\'etudier ici le cas de fibres quelconques.
\end{quote}
Des espaces \og fibr\'es\fg, qui, en termes modernes, seraient plut\^ot qualifi\'es de \og feuillet\'es\fg, sont en effet apparus dans la th\`ese de Seifert~\citepyear{Seifert33}. L'espace total y est une vari\'et\'e de dimension~$3$, les \og fibres\fg sont des cercles, la vari\'et\'e est localement diff\'eomorphe \`a $\RR^3=\RR^2\times\RR$, les fibres ou feuilles correspondant par ce diff\'eomorphisme aux $\lefacc a\rigacc\times\RR$. Il y a une projection sur l'espace des feuilles. 

Un fibr\'e localement trivial est un espace $E$ muni d'une projection $p$ sur un autre espace $B$ de telle sorte que $B$ est recouvert par des ouverts $U$ tels que $p^{-1}(U)$ soit hom\'eomorphe \`a un produit $U\times F$ (par un hom\'eomorphisme compatible avec les projections, et pr\'eservant donc les fibres).

\begin{figure}[htb]
\centerline{\includegraphics[scale=.8]{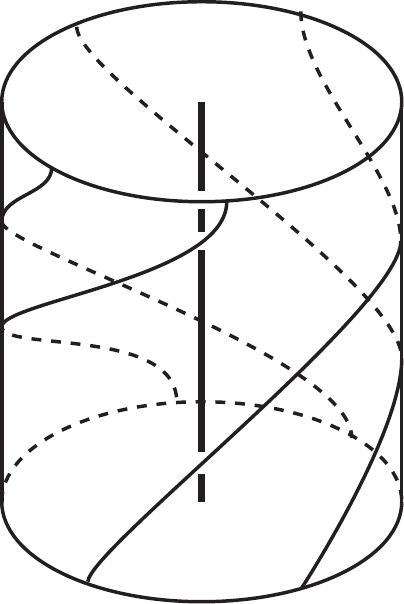}}
\caption{Un fibr\'e de Seifert}\label{figseifert}
\end{figure}

Dans le cas de Seifert, il n'y a pas de trivialit\'e locale, certaines fibres pouvant \^etre \og exceptionnelles\fg. Le mod\`ele local est celui du quotient du tore plein $D^2\times S^1$ par l'op\'eration du cercle
$$u\cdot(w,z)=(u^mw,u^nz), \text{ avec }m,n\in\ZZ\text{ premiers entre eux}.$$ 
Dans cette formule, $u$, $v$ et $z$ sont tous des nombres complexes, $u\in S^1$ (le cercle unit\'e des nombres complexes de module~$1$), de m\^eme que $z$, quand \`a $w$, c'est un point du disque unit\'e, un nombre complexe de module $\leq 1$. La fibre \og centrale\fg $w=0$ est une fibre exceptionnelle. Sur la figure \ref{figseifert}, $m=1$ et $n=2$; la fibre centrale est repr\'esent\'ee en gras, quelques autres fibres sont esquiss\'ees, le cylindre est bien de la forme $D^2\times [0,1]$ (les fibres forment un feuilletage), mais si on en recolle les deux extr\'emit\'es, on obtient un tore que l'on ne peut pas \'ecrire comme un produit, m\^eme localement.

Dans son article~\citepyear{Whitney37}, Whitney a consid\'er\'e, lui, des fibr\'es qui, comme ceux de Feldbau, sont \og localement triviaux\fg, mais dont les fibres sont des sph\`eres. Il a d'ailleurs donn\'e un \og important expos\'e sur les fibr\'es en sph\`eres\fg lors d'un congr\`es \`a Moscou en 1935 (dit Weil dans son livre~\citepyear[p.~115]{weilsouvenirs}). \`A ces r\'ef\'erences, il faut adjoindre celle aux travaux de de Rham, qu'\'Elie Cartan a fait ajouter \`a Feldbau\footnote{Le math\'ematicien suisse Georges de Rham (1903--1990) a pass\'e sa th\`ese \`a Paris en 1931. C'est la main d'\'Elie Cartan (qui connaissait bien de Rham) qui a ajout\'e cette note sur le manuscrit. \`A ma connaissance, de Rham n'a jamais publi\'e ces r\'esultats. Il ne les mentionne d'ailleurs pas dans son article de souvenirs de cette \'epoque~\citepyear{deRhamsouvenirs}.}:
\begin{quote}
J'ai appris que M.~de Rham a obtenu \`a peu pr\`es les m\^emes r\'esultats, qui \`a ma connaissance n'ont pas encore \'et\'e publi\'es.
\end{quote}
Tout ceci \'etait tr\`es r\'ecent. Remarquons encore que, dans cette note, les fibres sont des vari\'et\'es compactes (g\'en\'eralisant les cercles de Seifert et les sph\`eres de Whitney). Pour Feldbau comme pour Ehresmann, comme pour Seifert et Whitney, les espaces consid\'er\'es sont toujours des vari\'et\'es. \`A cette \'epoque, la distinction entre g\'eom\'etrie diff\'erentielle et topologie n'existe pas encore.

Les r\'esultats que d\'emontre Feldbau sont les suivants (en langage \`a peine modernis\'e).
\begin{itemize}
\item \emph{Un fibr\'e sur une base $B$ qui est un simplexe est trivialisable} (c'est le th\'eor\`eme A). Le fibr\'e est localement trivial, donc trivial sur d'assez petits simplexes. Le lemme-clef est alors le fait que, si le fibr\'e est trivial sur deux simplexes ayant une face commune, il est trivial sur leur r\'eunion. 
\item \emph{Les classes d'isomorphisme de fibr\'es sur la sph\`ere $S^n$ sont en bijection avec les classes d'homotopie d'applications de la sph\`ere $S^{n-1}$ dans le groupe $G$ des automorphismes de la fibre} --- au moins si ce groupe est connexe\footnote{Une hypoth\`ese que Feldbau omet de faire dans cette note, mais il corrigera dans~\cite{Laboureur42}.} --- (c'est le th\'eor\`eme B). Chacun des deux h\'{e}misph\`{e}res de la sph\`ere est hom\'eomorphe \`a un disque, donc \`a un simplexe. Le fibr\'e est donc trivialisable sur chacun des deux h\'emisph\`eres de $S^n$. Pour le reconstituer, il suffit de recoller les deux trivialisations le long de l'\'equateur $S^{n-1}$, soit de donner une application de $S^{n-1}$ dans $G$. Et, comme il sait que le fibr\'e tangent \`a la sph\`ere $S^{2n}$ n'est pas trivialisable, Feldbau en d\'eduit que le groupe $\pi_{2n-1}(\SO(2n))$ des classes d'homotopie d'applications de l'\'equateur $S^{2n-1}$ dans $\SO(2n)$ n'est pas trivial. \'Evaluer un groupe d'homotopie n'est pas toujours une mince affaire et on en connaissait assez peu \`a cette \'epoque.
\end{itemize}

\begin{figure}[htb]
\centerline{\hfill\includegraphics[scale=.8]{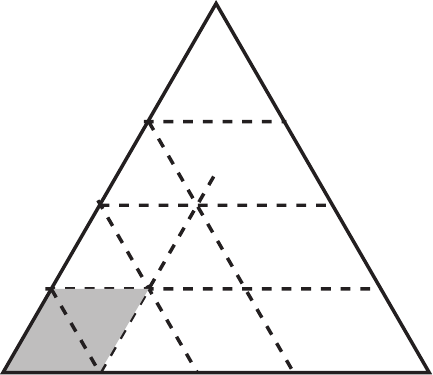}\hfill\includegraphics[scale=.8]{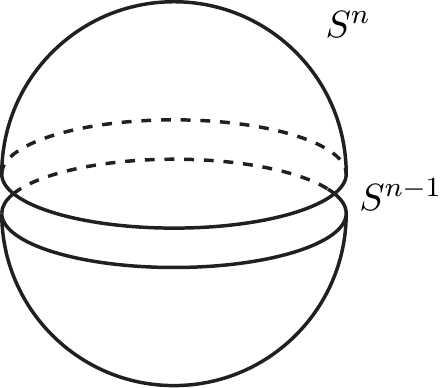}\hfill}
\caption{Le th\'eor\`eme A et le th\'eor\`eme B}
\end{figure}

\subsubsection{Sur les propri\'et\'es d'homotopie des espaces fibr\'es}
La deuxi\`eme note~\cite{ehresmannfeldbau} est \'ecrite par Feldbau en collaboration avec son directeur de th\`ese Ehresmann (elle est encore pr\'esent\'ee par \'Elie Cartan \`a l'Acad\'emie des sciences lors de la s\'eance du 4 juin 1941). Les auteurs y montrent des r\'esultats de base de la th\'eorie des fibr\'es: le rel\`evement des homotopies, ce que l'on formule aujourd'hui comme la suite exacte d'homotopie d'une fibration $p:E\to B$, 
$$\MRE{}\pi_n(F)\MRE{}\pi_n(E)\MRE{}\pi_n(B)\MRE{}\pi_{n-1}(F)\MRE{}$$
qui sont publi\'es simultan\'ement par Eckmann ou Hurewicz et Steenrod\footnote{Le math\'ematicien suisse Beno Eckmann est n\'e en 1917. Witold Hurewicz (1904--1956), n\'e en Pologne, a pass\'e sa th\`ese \`a Vienne en 1926, a \'emigr\'e aux \'Etats-Unis avant la guerre. C'est un des inventeurs des groupes d'homotopie \og sup\'erieurs\fg, c'est-\`a-dire des $\pi_n$ avec $n\geq 2$. Norman Steenrod (1910--1971) a pass\'e sa th\`ese \`a Princeton avec Lefschetz.} --- \og ce qui, compte tenu des circonstances pr\'esentes, \'etait inconnu des auteurs [Ehresmann et Feldbau]\fg, comme le dit \'el\'egamment Weil depuis l'Am\'erique dans sa recension pour \emph{Math. Reviews} en~1942.

\subsection*{Entre temps}Entre temps, c'est-\`a-dire dans l'intervalle de deux ans entre les deux notes~\cite{Feldbau39} et~\cite{ehresmannfeldbau}, il y a eu, bien s\^ur, la guerre. Feldbau l'a faite comme aviateur\footnote{Pour des informations sur la vie (et la mort) de Jacques Feldbau, voir~\cite{georgescerf,feldbau}.}. D\'emobilis\'e apr\`es l'armistice de juin 1940, il a \'et\'e nomm\'e professeur agr\'eg\'e au lyc\'ee de Ch\^ateauroux \`a la rentr\'ee 1940. Pas pour longtemps puisque le premier \og statut des juifs\fg, promulgu\'e par Vichy le 3 octobre (voir l'introduction de cet article), a entra\^in\'e sa r\'evocation de l'enseignement. Il a alors rejoint \`a Clermont-Ferrand l'universit\'e de Strasbourg repli\'ee et s'est remis au travail avec Ehresmann.

\subsubsection{Espaces fibr\'es associ\'es}\label{subsubespacesfibres}
Nous avons introduit parmi les publications de Jacques Feldbau, \`a la suite de Michel Zisman dans~\citepyear{zismanjames}, une note~\cite{Ehresmann41}, sign\'ee du seul Ehresmann, toujours pr\'esent\'ee par \'Elie Cartan, le 27 octobre~1941, mais dont la premi\`ere phrase affirme que 
\begin{quote}
Les r\'esultats qui vont \^etre expos\'es sont dus \`a la collaboration de l'auteur et d'un de ses \'el\`eves; ils font suite \`a une Note ant\'erieure.
\end{quote}
Le contenu d\'esignant la note ant\'erieure (c'est~\cite{ehresmannfeldbau}) et l'\'el\`eve en question, \'etait bien suffisant pour ajouter cette note \`a la liste de publications de Feldbau. 

\medskip
La consultation du manuscrit, conserv\'e dans la pochette de s\'eance du 1\up{er} d\'ecembre, a confirm\'e cette attribution. En r\'ealit\'e la note a \'et\'e envoy\'ee \`a \'Elie Cartan, comme la pr\'ec\'edente, avec les deux noms d'auteurs Charles Ehresmann et Jacques Feldbau. Les noms des auteurs ont \'et\'e ray\'es sur le manuscrit et remplac\'es, apparemment par la main d'\'Elie Cartan, par celui du seul Charles Ehresmann, la m\^eme main a d'abord ajout\'e 
\begin{quote}
Les r\'esultats qui vont \^etre expos\'es sont dus \`a la collaboration de l'auteur et de M.~Jacques Feldbau.
\end{quote}
puis a ray\'e cette mention et l'a remplac\'ee par celle qui a \'et\'e publi\'ee et dans laquelle le nom de Feldbau n'appara\^it plus: pour savoir qui est l'\'el\`eve en question, il faut aller lire la note ant\'erieure --- on ne fait pas plus prudent. Ce sont peut-\^etre ces h\'esitations\footnote{Ce n'\'etait certes pas une d\'ecision facile, surtout peut-\^etre pour \'Elie Cartan, que Camille Marbo (c'est-\`a-dire Marguerite Borel) a d\'ecrit dans son livre de souvenirs~\citepyear[p.~171]{marbo} comme un \og math\'ematicien transcendant, tr\`es courageux, tr\`es doux, mais d'un caract\`ere h\'esitant\fg.} avant la d\'ecision de publier la note, m\^eme sans le nom d'un de ses auteurs, qui sont responsables de l'inhabituel d\'elai entre la date de pr\'esentation de la note (le 27 octobre) et la date de la s\'eance dans laquelle elle est publi\'ee (le 1\up{er} d\'ecembre 1941).

Visiblement, m\^eme si \'Elie Cartan souhaitait publier cette note, l'Acad\'emie des sciences ne voulait pas de Feldbau dans ses publications. Et c'est ainsi qu'il dispara\^it comme auteur.

Venons-en au contenu. C'est dans cette note que les notions fondamentales de fibr\'e associ\'e et de fibr\'e principal sont d\'efinies. Le simplexe sur lequel tous les fibr\'es sont trivialisables depuis le th\'eor\`eme A de~\cite{Feldbau39} devient ici un espace contractile plus g\'en\'eral. Une occasion de mentionner que le nom de Feldbau aurait m\^eme pu dispara\^itre compl\`etement de cet article, comme on le voit sur le manuscrit justement quand l'auteur et son \'el\`eve signalent dans une note infrapaginale qu'ils g\'en\'eralisent un th\'eor\`eme de Jacques Feldbau. La partie de cette note comportant le nom interdit a \'et\'e ray\'ee elle aussi sur le manuscrit, pourtant la petite note figure compl\`etement, avec le nom de l'auteur du th\'eor\`eme A, dans le texte publi\'e. 

Dans la ligne de~\cite{ehresmannfeldbau,Ehresmann41}, Ehresmann publiera en 1942 une troisi\`eme note~\cite{Ehresmann42} qui \og emploie les d\'efinitions et notations de deux notes ant\'erieures\fg.

\subsubsection{Les structures fibr\'ees sur la sph\`ere et le probl\`eme du parall\'elisme} 
Nous en venons \`a l'article~\cite{Laboureur42}, le premier que Feldbau publie sous le pseudonyme transparent de Jacques Laboureur (Feldbau veut dire \og agriculture\fg, labourage, comme on disait autrefois, en allemand). Il s'agit d'une communication pr\'esent\'ee \`a la section de Clermont-Ferrand de la Soci\'et\'e math\'ematique de France le 16 avril 1942. Cette section se r\'eunit assez r\'eguli\`erement, Clermont est encore en zone \og libre\fg au printemps 42 (et jusqu'au 11 novembre)\footnote{Clermont-Ferrand est, pendant cette p\'eriode, le lieu d'une activit\'e math\'ematique incroyable. Voir par exemple~\cite{coutyglaeserperol,schwartz}.}. Ce printemps-l\`a, cette section discute, le 16 avril d'espaces fibr\'es (elle entend des communications de de Rham, d'Ehresmann et de \og Laboureur\fg), se r\'eunit \`a nouveau le 4 mai (Schwartz) et le 21 mai (Roussel, Delange). Dans cet article, Feldbau corrige l'\'enonc\'e du th\'eor\`eme B mentionn\'e ci-dessus (qui n'est vrai, comme je l'ai signal\'e, que si $G$ est connexe). Il y \'etudie aussi, en termes homotopiques, la parall\'elisabilit\'e de la sph\`ere~$S^n$ (c'est-\`a-dire la question de savoir si le fibr\'e tangent \`a la sph\`ere est trivialisable).

\subsubsection{Propri\'et\'es topologiques des automorphismes de la sph\`ere}
Le dernier article publi\'e du vivant de Jacques Feldbau\footnote{Arr\^et\'e au cours de la rafle des \'etudiants strasbourgeois \`a Clermont-Ferrand le 25 juin 1943, Jacques Feldbau sera envoy\'e \`a Drancy puis \`a Auschwitz-Monowitz en octobre 1943. Il survivra \`a la mortelle \'evacuation du  camp de janvier 1945, mais mourra d'\'epuisement juste avant la fin de la guerre au camp de Ganacker. Voir~\cite{georgescerf,feldbau}.} est~\cite{Laboureur43} dans lequel, toujours sous le nom de Laboureur, il \'etudie les relations entre le groupe des hom\'eomorphismes de degr\'e $+1$ de la sph\`ere dans elle-m\^eme et le groupe des rotations. Un hom\'eomorphisme de degr\'e $+1$ est un hom\'eomorphisme qui conserve l'orientation. Les rotations, c'est-\`a-dire les \'el\'ements du groupe $\SO(n+1)$ en font partie. Feldbau croyait avoir d\'emontr\'e que le groupe $\homeo(S^n)$ se r\'etracte par d\'eformations sur son sous-groupe, le groupe orthogonal $\On(n+1)$. Le recenseur de \emph{Math. Reviews}, Hans~Samelson\footnote{Le topologue Hans Samelson (1916--2005) est n\'e \`a Strasbourg comme Feldbau, mais lui \'etait allemand, il a fait ses \'etudes \`a Breslau puis \`a Z\"urich (avec Hopf), et a \'emigr\'e aux \'Etats-Unis en 1941.}, l'avait d\'ej\`a remarqu\'e, il y a une erreur dans la d\'emonstration; on sait d'ailleurs aujourd'hui que ce r\'esultat est faux pour $n\geq 5$ (voir, pour une \'etude plus pr\'ecise mais lisible par des non-sp\'ecialistes, l'article de Douady~\citepyear{douadycerf} sur les travaux de Jean Cerf, dans lequel l'\'enonc\'e de~\cite{Laboureur43} est nomm\'e \og conjecture de Feldbau\fg).

\subsubsection{Sur la loi de composition entre \'el\'ements des groupes d'homotopie} Le dernier article~\cite{Feldbau58} est tardif (et posthume). On pourra consid\'erer qu'il n'entre pas vraiment dans le cadre du pr\'esent travail puisqu'il est paru en 1958. Il s'agit quand m\^eme d'un travail datant de la p\'eriode envisag\'ee ici et qui n'a pas pu \^etre publi\'e en son temps (parce que son auteur avait \'et\'e d\'eport\'e). Feldbau y \'etudie les applications d'un produit de sph\`eres dans un espace topologique. Il d\'efinit d'abord ce que l'on appelle aujourd'hui le \og produit de Whitehead\fg, puisque Whitehead\footnote{Le math\'ematicien anglais John Henry Constantine Whitehead (1904--1960) a fait sa th\`ese \`a Princeton en 1932 avec Veblen.} l'a d\'efini lui aussi et a pu publier ce travail~\citepyear{WhiteheadProduit} pendant la guerre,
$$\pi_p(X)\times\pi_q(X)\MRE{}\pi_{p+q-1}(X)$$
($X$ est un espace topologique point\'e, dont le point base sera not\'e $\star$ dans la suite). Il \'ecrit un \og diagramme de Heegaard\fg pour la sph\`ere $S^{p+q-1}$, c'est-\`a-dire qu'il repr\'esente cette sph\`ere comme
\begin{figure}[htb]
\centerline{\includegraphics[scale=.8]{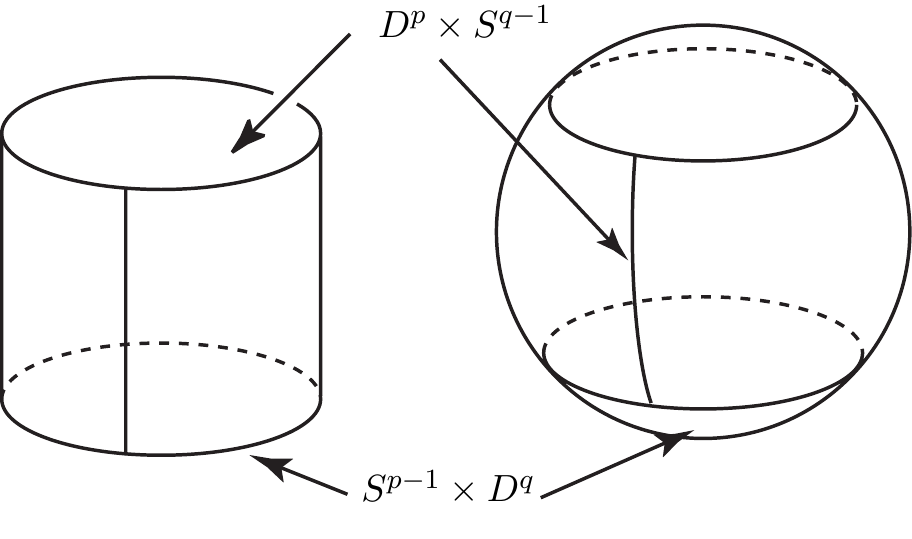}}
\caption{Le produit de Whitehead}
\end{figure}

$$S^{p+q-1}=S^{p-1}\times D^q\cup_{S^{p-1}\times S^{q-1}}D^p\times S^{q-1}=\partial(D^p\times D^q).$$
Dans le cas de la sph\`ere $S^3$ (c'est-\`a-dire si $p=q=2$), il s'agit du diagramme de Heegaard standard de genre~$1$ qui repr\'esente cette sph\`ere comme r\'eunion de deux tores pleins. La figure, elle, repr\'esente le cas de la sph\`ere $S^2$ (avec $p=1$ et $q=2$). 

Si les applications 
$$f:(D^p,S^{p-1})\MRE{}(X,\star)\text{ et }g:(D^q,S^{q-1})\MRE{}(X,\star)$$
repr\'esentent un \'el\'ement $(\alpha,\beta)$ de $\pi_p(X)\times\pi_q(X)$, il d\'efinit l'image $[\alpha,\beta]$ comme la classe d'homotopie de
\begin{align*}
h:S^{p+q-1}&\MRE{}X\\
z&\DMRE{}\begin{cases}
g(y)\text{ si }z=(x,y)\in S^{p-1}\times D^q\\
f(x)\text{ si }z=(x,y)\in D^p\times S^{q-1}
\end{cases}
\end{align*}
(la notation vient du fait que, si $p=q=1$, $[\alpha,\beta]$ est le commutateur $\alpha\beta\alpha^{-1}\beta^{-1}\in\pi_1(X)$). Le cas o\`u $p=1$ redonne l'op\'eration classique du groupe fondamental sur les groupes d'homotopie. 

Il applique ensuite cette construction \`a l'\'etude des applications d'un produit de sph\`eres dans l'espace topologique $X$. Dans la toute derni\`ere partie, il utilise le cas $p=q=n$, pour \'etudier $\pi_{2n-1}(S^n)$... en utilisant un r\'esultat annonc\'e par Freudenthal\footnote{Le topologue Hans Freudenthal (1905--1990), a fait une th\`ese avec Hopf \`a Berlin en 1931, puis toute sa carri\`ere aux Pays-Bas.} en 1939, \`a savoir, pour tout $n$ pair, l'existence d'une application $S^{2n-1}\to S^n$ d'invariant de Hopf $1$ (ce qui est faux sauf pour $n=2,4,8$). La toute derni\`ere partie de ce travail inachev\'e s'intitule \og Groupes de Ker\'ekj\'art\'o\fg, il s'agit  des groupes d'homotopie toro\"idaux (inspir\'es par un passage du livre~\cite{kerekjarto}).

\subsection*{Et apr\`es?}
Les fibr\'es sur les espaces contractiles sont toujours trivialisables et les topologues font si grand usage de cette propri\'et\'e que peu d'entre eux savent qu'il s'agit d'un th\'eor\`eme de Feldbau. Ils utilisent les fibr\'es associ\'es comme ils respirent et, dans le meilleur des cas, pensent savoir que c'est Ehresmann qui les a invent\'es.

C'est sans doute en partie gr\^ace \`a l'intervention de certains de ses amis que l'article posthume~\cite{Feldbau58} de Feldbau a finalement \'et\'e publi\'e (voir~\cite{feldbau}). Ehresmann l'a fait pr\'ec\'eder d'une pr\'eface; il a aussi \'ecrit, \`a la m\^eme \'epoque en 1958, un rapport sur les travaux de Feldbau (qui est reproduit dans le premier volume~\cite{EhresmannOC1} de ses \OE uvres compl\`etes) dans lequel il a fait la liste de ses articles... sans mentionner dans aucun de ces deux textes, pour une raison que je n'arrive pas \`a imaginer, que Feldbau est un des auteurs de la note~\cite{Ehresmann41}.

\medskip
On se reportera, pour l'histoire de la topologie et des fibr\'es, au livre~\cite{dieudonne}, pour l'homotopie \`a l'article~\cite{eckmann}, et surtout au livre~\cite{jameshistory} (notamment \`a l'article~\citepyear{zismanjames} de Michel Zisman que ce dernier livre contient et que j'ai abondamment utilis\'e pour \'ecrire ce paragraphe).

\section{L'Acad\'emie des sciences}\label{secacad}

Je n'ai pas su trouver, \`a l'Acad\'emie des sciences, d'expression conserv\'ee de la mani\`ere dont les directives officielles, celles du \textsc{mbf} (cit\'ees dans l'introduction de cet article), de Fourneau ou d'autres, ont \'et\'e transmises, discut\'ees, accept\'ees, ni dans ce qui para\^it, publiquement (dans les \emph{Comptes rendus}) ni dans ce qui est officiellement conserv\'e (les proc\`es-verbaux des r\'eunions non publiques des acad\'emiciens, dites \og en comit\'e secret\fg). Il n'y a rien eu d'explicite dont il reste des traces \'ecrites dans les Archives de l'Acad\'emie des sciences (voir la note~\ref{notedestruction}). 

\medskip
Une trace t\'enue conserv\'ee par l'Acad\'emie des sciences appara\^it sur une carte envoy\'ee par les \'editions Masson (qui publient les travaux de l'Acad\'emie de m\'edecine) en mars 1942 et conserv\'ee dans la pochette de la s\'eance du 23 mars 1942. En voici le texte:
\begin{quote}
Par d\'ecision du \emph{Milit\"arbefehlshaber in Frankreich}, qui nous a \'et\'e notifi\'ee par le Groupement Corporatif de la Presse P\'eriodique, le r\'egime de publication des p\'eriodiques m\'edicaux et scientifiques est profond\'ement r\'eduit et modifi\'e.

Les instructions que nous avons re\c cues ne sont pas encore d\'efinitives, nous ne pouvons donc pas encore mettre au point les mesures impos\'ees par la situation nouvelle, ni par cons\'equent les porter \`a la connaissance de nos abonn\'es.

Nous vous demandons de nous faire confiance dans la p\'eriode difficile que nous traversons.

\begin{flushright}
MASSON \&\ Cie
\end{flushright}
\end{quote}

Serait-ce la seule trace \'ecrite dans ces archives d'une d\'ecision sur les publications scientifiques prise par le \textsc{mbf}, et qui n'arriverait \`a l'Acad\'emie que par l'interm\'ediaire du Groupement Corporatif puis des \'editions Masson?

\subsection*{Les difficult\'es mat\'erielles}
Je ne mentionne pas ici les difficult\'es de la vie quotidienne des chercheurs, mais celles li\'ees \`a la publication de leurs travaux. Les communications sont difficiles. Par exemple, beaucoup de scientifiques sont en zone \og libre\fg, notamment de nombreux math\'ematiciens \`a Clermont-Ferrand. L'Acad\'emie des sciences est, elle, \`a Paris, on communique donc par cartes interzones et/ou par l'interm\'ediaire du Secr\'etaire d'\'Etat \`a l'\'Education Nationale et \`a la Jeunesse, \`a Vichy, qui transmet. Il y a m\^eme un exemple d'une note dont le manuscrit est arriv\'e sous la forme de cartes interzones num\'erot\'ees de 1 \`a 6 (note~\citepyear{Mirguet} de g\'eom\'etrie infinit\'esimale, par Jean Mirguet, pochette de la s\'eance du 4 ao\^ut 1941). Plusieurs notes re\c cues par l'Acad\'emie des sciences ont m\^eme \'et\'e exp\'edi\'ees de camps de prisonniers de guerre (celle de Fr\'ed\'eric Roger\footnote{Il y aura plus tard, en 1942, trois notes de Jean Leray en provenance du m\^eme \emph{Oflag} XVII~A. Je n'en dis pas plus, comptant revenir ailleurs sur cette question (voir le paragraphe de conclusion de cet article).}, pr\'esent\'ee par Borel \`a la s\'eance du 13 janvier 1941, par exemple).

\subsection*{Les d\'ebats \`a l'Acad\'emie des sciences}
Si l'institution a choisi l'\og accommodement\footnote{Voir la note~\ref{accommodement}.}\fg, il est certain qu'il y a eu, entre ses membres, des d\'ebats, des discussions, des d\'esaccords politiques, m\^eme si aucun compte-rendu \'ecrit n'en a \'et\'e conserv\'e. J'ai d'abord \'et\'e extr\^emement \'etonn\'ee qu'il n'ait subsist\'e aucune trace dans les proc\`es-verbaux d'une quelconque r\'eaction de ses coll\`egues \`a l'arrestation de Paul Langevin\footnote{Le c\'el\`ebre physicien Paul Langevin (1872--1946) a \'et\'e, \`a cause de ses opinions anti-fascistes, incarc\'er\'e le 30 octobre 1940 par la Gestapo \`a la prison de la Sant\'e, lib\'er\'e quarante jours plus tard et assign\'e \`a r\'esidence \`a Troyes.}, un des membres les plus c\'el\`ebres de l'Acad\'emie des sciences, le 30 octobre 1940, ni \`a la r\'epression brutale de la manifestation des \'etudiants \`a l'Arc de Triomphe quelques jours plus tard, le 11 novembre 1940, qui a \'et\'e jusqu'\`a entra\^iner la fermeture de la Sorbonne. Il y a pourtant eu d\'ebat, il y a pourtant eu indignation. Ces \'ev\'enements ont boulevers\'e (au moins certains) des acad\'emiciens. Le math\'ematicien allemand Helmut Hasse, en visite dans le Paris occup\'e comme commis-voyageur de la science allemande\footnote{Helmut Hasse (1898--1979), directeur de l'Institut de G\"ottingen depuis 1934, \'etait le tr\'esorier de la \textsc{dmv} (\emph{Deutsche Mathematiker-Vereinigung}), il a fait pendant l'Occupation plusieurs visites \`a Paris dans son uniforme de capitaine de corvette (son grade dans la marine), notamment pour chercher des collaborateurs pour la revue de r\'ef\'erences \emph{Zentralblatt f\"ur Mathematik}. Pour les plans nazis de r\'eorganisation des math\'ematiques, voir~\cite{SiegmundSNeuordnung}.}, rapporte sa visite \`a \'Elie Cartan (Henri Cartan est pr\'esent) dans une lettre \`a S\"uss\footnote{\emph{Nachla\ss\ S\"u\ss}, universit\'e de Freiburg. Wilhelm S\"uss (1895--1958) \'etait le pr\'esident de la \textsc{dmv}, il sera le fondateur du \emph{Reichsinstitut} d'Oberwolfach. Voir par exemple~\cite{VRemmert}.} du 19 d\'ecembre 1940. L'atmosph\`ere \'etait plut\^ot r\'eserv\'ee,
\begin{quote}
\noindent\emph{Hier war die Stimmung etwas reserviert}
\end{quote} 
\`a cause principalement, pense-t-il, de deux \'ev\'enements politiques, la r\'epression de la manifestation des \'etudiants sur la tombe du soldat inconnu et l'arrestation de Langevin. Ces deux choses ont beaucoup \'emu Cartan \begin{quote}
\noindent\emph{Beide Dinge regten Cartan ersichtlich sehr auf}.
\end{quote} 
Pourtant, il a \'et\'e re\c cu gentiment
\begin{quote}
\noindent\emph{Pers\"onlich war er aber sehr nett}.
\end{quote} 
Dans une lettre \`a Krull\footnote{\emph{Nachla\ss\ Hasse}, universit\'{e} de G\"ottingen. L'alg\'{e}briste allemand Wolfgang Krull (1899--1971) \'etait un ami et coll\`egue de Hasse.} \'ecrite quelques semaines plus tard (le 23 janvier 1941), Hasse rapporte la m\^eme visite \`a Paris en disant que chez les Cartan l'atmosph\`ere \'etait d\'eprim\'ee et tendue
\begin{quote}
\noindent\emph{Die Stimmung bei beiden war etwas gedr\"uckt und auch ge\-spannt}
\end{quote} 
--- ses avances n'ont pas eu \'enorm\'ement de succ\`es, puisqu'il \'ecrit un peu plus loin que Julia\footnote{L'acad\'emicien Gaston Julia (1893--1978) a en effet cherch\'e et trouv\'e des collaborateurs fran\c cais pour \emph{Zentralblatt}, comme la suite de la correspondance de Hasse le montre.} est sans doute le seul math\'ematicien fran\c cais qui fera tout son possible pour la reprise des relations scientifiques avec l'Allemagne.

\medskip
Je me suis donc faite \`a l'id\'ee de ce silence (assourdissant) et me suis consacr\'ee \`a la recherche de traces des brisures de ce silence.

Il suffit de compter les r\'eunions en \og comit\'e secret\fg et leurs dur\'ees (les \emph{Comptes rendus} donnent l'heure de d\'ebut et l'heure de fin) et de confronter ces dur\'ees \`a la bri\`evet\'e des proc\`es-verbaux report\'es sur le registre des comit\'es secrets pour \^etre assur\'e que nos acad\'emiciens ont bien discut\'e. Il y a donc eu de nombreuses r\'eunions dont les d\'ebats sont rest\'es vraiment secrets. Le registre des comit\'es secrets fait \'etat de sept r\'eunions pour toute l'ann\'ee 1940, de six pour 1941... mais de vingt-six pour 1942 et de vingt pour 1943. Au-del\`a des informations sur les trois st\`eres de bois de chauffage en provenance de Chantilly auxquelles ont droit les acad\'emiciens (18 mai 1942) et sur l'\'electricit\'e suppl\'ementaire qui ne peut leur \^{e}tre accord\'ee \`a cause de la p\'enurie (15 novembre 1943), il y a forc\'ement eu des raisons s\'erieuses \`a la tenue de toutes ces r\'eunions, et \`a leur tenue en comit\'e secret (c'est-\`a-dire non public)\footnote{Pour tout ce paragraphe, registre des comit\'es secrets, archives de l'Acad\'emie des sciences.}.

Le dossier aux archives de l'Acad\'emie des sciences du math\'ematicien et acad\'emicien Paul Montel\footnote{Paul Montel (1876--1975), le math\'ematicien sp\'ecialiste d'analyse complexe qui a invent\'e les familles normales, a \'et\'e doyen de la facult\'e des sciences de Paris de 1941 \`a 1946.}, qui \'etait aussi doyen de la facult\'e des sciences, contient plusieurs cartes qu'il a envoy\'ees au Secr\'etaire Perp\'etuel de l'Acad\'emie des sciences l'informant au fur et \`a mesure des arrestations de coll\`egues scientifiques, acad\'emiciens ou pas (arrestations de Mauguin, Langevin, Borel en octobre 1941, arrestations de Caullery, Cotton, Dussaud, Lacaut, Ferdinand Lot, Pelliot, Petit-Dutaillis, et de Marcel Bouteron, biblioth\'ecaire de l'Institut et Georges Bourgin, archiviste). Il est certain que ces messieurs en ont parl\'e.

Voici quelques dates de ces d\'ebats\footnote{Pour cette \'enum\'eration, registre des comit\'es secrets et pochettes des s\'eances correspondantes, archives de l'Acad\'emie des sciences.}.
\begin{itemize}
\item Le 10 juin 1940 (juste avant l'armistice, donc). La discussion porte sur la question de savoir si l'Acad\'emie restera \`a Paris ou non en cas d'occupation, il y a tr\`es peu d'acad\'emiciens pr\'esents (dix-huit votants) et aucune d\'ecision n'est prise.
\item Le 19 ao\^ut 1940. Les acad\'emiciens d\'ebattent des rapports avec la censure allemande, l'occasion est le discours prononc\'e par le Pr\'esident, le g\'en\'eral Perrier\footnote{Le pr\'esident Georges Perrier (1872--1946) \'etait un g\'eographe.}, \`a l'occasion de la mort de Breton\footnote{Jules Breton (1872--1940), ing\'enieur chimiste, acad\'emicien libre depuis 1920, directeur des inventions int\'eressants la D\'efense Nationale pendant la 1\up{\`ere} guerre mondiale (mise au point du char d'assaut), ministre de l'hygi\`ene en 1920.}, et que les acad\'emiciens ont pr\'ef\'er\'e amputer de quelques phrases sur le r\^ole jou\'e par Breton \og pendant la guerre pr\'ec\'edente\fg, une des raisons invoqu\'ees \'etant de ne pas attirer l'attention des Allemands sur les installations de Bellevue et sur ce qui se fait dans certains laboratoires du \textsc{cnrs}. Ceci a entra\^in\'e la d\'emission de Perrier de la pr\'esidence et ne pouvait donc pas \^etre cach\'e. Ce qui reste de ces d\'ebats: une version manuscrite (prise de notes pendant la r\'eunion), deux versions dactylographi\'ees dont celle qui a \'et\'e conserv\'ee dans le registre des comit\'es secrets, un peu \'edulcor\'ee (la censure redout\'ee n'est plus qualifi\'ee d'allemande, les interventions des uns et des autres ne sont pas d\'etaill\'ees).

\item Le 4 novembre 1940. Il s'agit d'une discussion sur la mani\`ere de transmettre les \emph{Comptes rendus} en zone \og libre\fg et \`a l'\'etranger. Les notes manuscrites conserv\'ees dans la pochette de la s\'eance, quoiqu'extr\^emement concises, montrent une discussion plus riche que ce qui appara\^it dans le registre des comit\'es secrets: il n'est pas totalement innocent de proposer de passer par l'Ambassade des \'Etats-Unis (Aim\'e Cotton), celle d'Espagne (\'Emile Borel, s\'eparant la question de l'\'etranger de celle de la zone libre), par le Minist\`ere (\'Emile Picard)... jusqu'\`a: \og On pourrait demander aux Allemands eux-m\^emes (Deslandres). --- Nous ne pouvons avoir de rapport direct avec les occupants\fg (Picard)\footnote{Le physicien Aim\'e Cotton (1869--1951) \'etait un homme de gauche, il sera d\'ecor\'e de la m\'edaille de la R\'esistance. Il en est de m\^eme du math\'ematicien \'Emile Borel (1871--1956), qui a aussi \'et\'e d\'eput\'e (radical) de 1924 \`a 1936 et ministre en 1925. \'Emile Picard (1856--1941) \'etait Secr\'etaire Perp\'etuel depuis 1917. L'astronome Henri Deslandres (1853--1948) avait dirig\'e l'Observatoire de Paris de 1927 \`a 1929.}.

\item 
Il y a eu des affrontements politiques, comme celui qui a vu, d\'ebut 1942, apr\`es la mort du Secr\'etaire perp\'etuel \'Emile Picard, le retrait sans doute impos\'e par le pouvoir de l'ind\'esirable \'Emile Borel\footnote{Dossier biographique d'\'Emile Borel, archives de l'Acad\'emie des sciences.} et l'\'election du plus acceptable Louis de Broglie\footnote{Le physicien Louis de Broglie (1892--1987), qui a re\c cu le Prix Nobel de physique en 1929 pour sa d\'ecouverte de la nature ondulatoire de l'\'electron, a \'et\'e \'elu \`a l'Acad\'emie des sciences en 1933.} --- un affrontement dont on pourrait tr\`es bien ne pas s'apercevoir. 

Picard meurt le 12 d\'ecembre 1941. Il faut \'elire un nouveau Secr\'etaire perp\'etuel. \'Emile Borel a \'et\'e arr\^et\'e le 10 octobre et incarc\'er\'e \`a Fresnes o\`u il a pass\'e cinq semaines dans des conditions tr\`es difficiles (rappelons qu'il avait alors soixante-dix ans), voir~\cite{marbo}. Le 12 janvier, Borel \'ecrit une lettre dans laquelle il explique (ou plus exactement d\'eclare) qu'il retire sa candidature \`a la succession de Picard \og \`a la suite d'une nouvelle discussion avec M. Alfred Lacroix\footnote{Le min\'eralogiste Alfred Lacroix  (1863--1948) \'etait l'\og autre\fg secr\'etaire perp\'etuel.}\fg. \`A l'anciennet\'{e}, les candidats naturels auraient \'et\'e Hadamard, \'elu en 1912 mais r\'efugi\'e aux \'Etats-Unis, Borel, \'elu en 1921, Cartan, \'elu en 1931, de Broglie \'elu en 1933. Les candidats sont finalement \'Elie Cartan et Louis de Broglie, et c'est ce dernier qui est \'elu. \`A la Lib\'eration, Borel essaiera d'obtenir que celui-ci d\'emissionne pour lui laisser la place qui lui revenait, \'eclairant ce qui s'est pass\'e en 1942 et qui aurait bien pu passer inaper\c cu\footnote{Il est remarquable aussi qu'il n'y ait aucune trace de cette non-\'election de Borel comme Secr\'etaire perp\'etuel, ni dans le livre de souvenirs de Camille Marbo (son \'epouse)~\citepyear{marbo}, ni dans la biographie~\cite{guiraldenq}.}:

\medskip
\begin{quote}
[...]
D\`es la lib\'eration, j'ai cherch\'e \`a \^etre exactement renseign\'e sur les r\'eactions de l'Acad\'emie au moment de mon incarc\'eration \`a Fresnes. M. Vincent\footnote{Hyacinthe Vincent (1862--1950) \'etait le pr\'esident de l'Acad\'emie des sciences en 1941.} a bien voulu me raconter en d\'{e}tail ses d\'emarches, en tant que Pr\'esident de l'Acad\'emie, aupr\`es de M. de Brinon\footnote{Fernand de Brinon (1885--1947) repr\'esente le gouvernement de Vichy aupr\`es du haut-commandement allemand \`a Paris.}. Celui-ci le persuada que toute d\'emarche en faveur des acad\'emiciens incarc\'er\'es ne pourrait avoir que les plus grands dangers pour eux-m\^emes et pour l'Acad\'emie. Par suite, lorsque, apr\`es le d\'ec\`es de M. Picard un courant se manifesta dans l'Acad\'emie en faveur de ma candidature, M. Vincent ainsi que M. Esclangon\footnote{Ernest Esclangon (1876--1954), directeur de l'Observatoire de Paris et vice-pr\'esident de l'Acad\'emie des sciences, sera pr\'esident en 1942.}, vice-pr\'esident, consid\'er\`erent comme de leur devoir de s'opposer \`a cette candidature. M. Vincent a r\'ep\'et\'e cette conversation \`a M. Roussy\footnote{Gustave Roussy (1874--1948), acad\'emicien, \'etait recteur de l'Acad\'emie de Paris jusqu'\`a sa r\'evocation apr\`es la manifestation \'etudiante du 11 novembre 1940.} et M. Esclangon me l'a confirm\'ee, en ajoutant, qu'\`a son avis, j'aurais \'et\'e s\^urement \'elu, si je n'avais pas \'et\'e arr\^et\'e par les Allemands. De l'avis de M. Vincent et de M. Esclangon, c'est donc cette arrestation seule qui a emp\^ech\'e mon \'election puisque tout candidat qui m'\'etait oppos\'e devait avoir, outre ses voix personnelles, toutes celles de ceux qui partageaient l'opinion des Pr\'esidents sur les dangers de mon \'election. Malgr\'e cela, au d\'ebut de janvier, j'avais de nombreuses promesses et mon succ\`es paraissait assur\'e. Le coup de gr\^ace me fut alors donn\'e par Monsieur Carcopino\footnote{J\'er\^ome Carcopino (1881--1970), assure les fonctions de recteur de l'Acad\'emie de Paris apr\`es la r\'evocation de Roussy, secr\'etaire d'\'Etat \`a l'\'Education nationale et \`a la jeunesse dans le gouvernement Darlan en 1942, directeur de l'\'Ecole Normale Sup\'erieure.} qui, comme vous le savez mieux que personne, exigea le retrait de ma candidature, pour les m\^emes raisons sugg\'er\'ees \`a M. Vincent par M. de Brinon.
 
Je suis donc en droit de penser que, sans le vouloir, M.~[illisible] m'a impos\'e une peine suppl\'ementaire, s'ajoutant aux cinq semaines d'incarc\'eration. Il me semble que j'ai le droit de demander qu'une r\'eparation me soit accord\'ee. La r\'eparation la plus compl\`ete serait la d\'emission de M. Louis de Broglie et mon \'election. 

[...]\footnote{Lettre de Borel au Secr\'etaire Perp\'etuel Alfred Lacroix, le 23 septembre 1944. Dossier biographique de Borel aux archives de l'Acad\'emie des sciences.}
\end{quote}

\item
D'autres sont sous-jacents au printemps 1942 lorsque la question de la reprise ou non de l'\'election de membres fait l'objet de plusieurs d\'ebats. Un tableau, \`a vrai dire assez confus, des acad\'emiciens pr\'esents \`a Paris est dress\'e: certains ont \og \'emigr\'e\fg (c'est le cas d'Hadamard), d'autres sont en \og zone libre\fg, l'un est prisonnier de guerre, lit-on aussi, mais il n'appara\^it pas dans le tableau. Le d\'ebat fait appara\^itre des points de vue oppos\'es:
\begin{quote}
M.~Borel explique le sens de la d\'ecision prise par la section de g\'eom\'etrie. \`A son avis, il appartient \`a l'Acad\'emie de prendre une d\'ecision g\'en\'erale. La section de g\'eom\'etrie a entendu simplement d\'eclarer que si cette d\'ecision \'etait positive, elle \'etait pr\^ete \`a pr\'esenter une liste de candidats, mais il ajoute que sur la question de la reprise, la section se trouvait partag\'ee, M. Julia et lui \'etant partisans de la diff\'erer, MM. Cartan et Montel \'etant d'avis d'\'elire\footnote{Finalement, la section de g\'eom\'etrie \'elira Denjoy pour remplacer Lebesgue, d\'ec\'ed\'e en 1941.}.
\end{quote}
Les motivations des positions prises par les uns ou les autres peuvent \^etre vari\'ees. \'Evoquant la question de l'\'election de nouveaux membres en remplacement des \'emigr\'es, des juifs ou des prisonniers, mais dans le cas du Coll\`ege de France, Philippe Burrin~\citepyear[p.~314]{burrin} fait \'etat de la position de P\'etain, qui est pour la suspension des remplacements, oppos\'ee \`a celle de l'administrateur Faral\footnote{Edmond Faral (1882--1958), sp\'ecialiste de la litt\'erature du Moyen \^Age, a \'et\'e administrateur du Coll\`ege de France de 1937 \`a 1955.}, qui est pour l'\'election de nouveaux professeurs, pour le maintien de la \og force fran\c caise\fg. Dans le d\'ebat \`a l'Acad\'emie des sciences, on voit ici Borel et Julia, qui ne sont certainement pas des amis politiques, adopter la m\^eme position. Mais revenons \`a la discussion, qui comporte une remarque sibylline:

\medskip
\begin{quote}
M.~Jacob\footnote{L'acad\'emicien g\'eologue Charles Jacob (1878--1962), a dirig\'e le \textsc{cnrs} de 1940 \`a 1944.} pense que l'Acad\'emie doit adopter une r\`egle g\'en\'erale et appelle l'attention sur le fait qu'il serait actuellement impossible de voter pour certains candidats m\'eritant cependant d'\^etre pr\'esent\'es en premi\`ere ligne.
\end{quote}

\medskip
Il va sans dire qu'il ne subsiste aucune trace \'ecrite de l'intervention d'un acad\'emicien collaborateur qui aurait demand\'e l'exclusion d'Hadamard sous pr\'etexte que celui-ci n'assistait pas aux s\'eances. On trouvera cette \og information\fg, que la tradition orale a fait conna\^itre \`a Schwartz, dans~\citepyear[p.~152]{schwartz} (Schwartz ne donne pas le nom de l'acad\'emicien en question).

\item 
Il y a eu un seul moment o\`u les opposants ont r\'eussi \`a faire appara\^itre une trace \'ecrite de l'application des lois antis\'emites, celui de l'\'election de Joliot, beaucoup plus tard. La section de physique est \`a l'\'epoque la seule section de l'Acad\'emie des sciences \`a \^etre \og \`a gauche\fg (c'est celle notamment d'Aim\'e Cotton, de Charles Fabry\footnote{Sur Charles Fabry (1867--1945) et les publications, voir ci-dessous le~\T\ref{subsecpollaczek}.} et de Paul Langevin). Non contente de faire \'elire Joliot \`a l'Acad\'emie des sciences (\`a la place laiss\'ee libre par la mort de Branly\footnote{Le physicien Fr\'ed\'eric Joliot-Curie (1900--1958), inventeur de la radioactivit\'e artificielle (et prix Nobel) avec Ir\`ene Joliot-Curie, va ainsi \^etre \'elu... au fauteuil pour lequel l'Acad\'emie des sciences avait cru bon de pr\'ef\'erer Branly \`a Marie Curie.}), la section de physique\footnote{Michel Pinault m'a dit penser qu'il s'agissait d'une provocation de la section de physique g\'en\'erale pour forcer \`a un d\'ebat sur ce th\`eme. Pour une passionnante \'etude des circonstances de cette \'election et de son contexte, voir son livre~\citepyear{pinaultlivre}, en particulier les pages~219 et suivantes.} a r\'eussi \`a faire en sorte qu'apparaissent, dans le proc\`es-verbal du Comit\'e Secret qui a pr\'epar\'e cette \'election, le 21 juin 1943,
\begin{itemize}
\item outre les noms que l'on peut lire dans les \emph{Comptes rendus} (Joliot en premi\`ere ligne, en deuxi\`eme ligne Becquerel, Cabannes, Ribaud), ceux d'Henri Abraham et d'Eug\`ene Bloch,
\item et donc la mention \og M. A. Lacroix fait observer que cette liste comprend les noms de deux savants dont les travaux sont fort estim\'es, mais dont les lois actuellement en vigueur ne permettent pas l'\'election. Il pense que ces noms ne peuvent sans inconv\'enients, peut-\^etre graves\footnote{C'est \`a des inconv\'enients pour l'Acad\'emie que pense, sans doute, Lacroix. Ce sont des inconv\'{e}nients graves qui vont frapper les physiciens Henri Abraham (1868--1943) et Eug\`ene Bloch (1878--1944), qui mourront tous deux en d\'eportation, le premier d\`es cette m\^eme ann\'ee, \`a soixante-quinze ans, le deuxi\`eme l'ann\'ee suivante, \`a soixante-six ans.}, \^etre maintenus sur la liste.\fg
\end{itemize}
\`A ma connaissance, c'est la seule occasion o\`u \og les lois actuellement en vigueur\fg sont explicitement mentionn\'ees, presque trois ans apr\`es la promulgation du premier \og statut des juifs\fg, mais pas exactement en relation avec les publications.

\end{itemize}

\subsection*{Les traces}
Il reste quand m\^eme des indices de l'application de ces consignes, d'abord dans les listes des auteurs des notes elles-m\^emes mais surtout dans ce qui a \'et\'e conserv\'e de telle ou telle s\'eance. C'est ce que j'essaie de montrer maintenant.

Au fil de l'ann\'ee 1941, on voit cette interdiction aux juifs de publier s'appliquer de plus en plus s\'ev\`erement, de plus en plus strictement. En voici quelques traces (chronologiquement). Le cas exceptionnel d'Andr\'e Bloch sera consid\'er\'e plus sp\'ecifiquement au~\T\ref{subsecbloch}. En mars paraissent deux notes de Laurent Schwartz~\citepyear{Schwartz41a, Schwartz41c} (les 3 et 17 mars) et une de Marie-H\'el\`ene Schwartz\footnote{C'est par erreur que \emph{Math. Reviews} attribue \`a Laurent Schwartz la note~\citepyear{MHSchwartz41} de Marie-H\'el\`ene Schwartz. La math\'ematicienne Marie-H\'el\`ene Schwartz, fille de Paul L\'evy et \'epouse de Laurent Schwartz, publie cette note sous le nom de \og M\up{me} Laurent Schwartz\fg. La table des mati\`eres du volume~{\bfseries 212} fait appara\^itre cette math\'ematicienne comme \textsc{Schwartz (M\up{me} Laurent), n\'ee Marie-H\'el\`ene L\'evy}, alors que, quelques semaines plus tard, les noms rep\'erables comme juifs vont \^etre masqu\'es et/ou supprim\'es.}~\citepyear{MHSchwartz41} (le 10 mars), ces trois notes sont transmises par Paul Montel; nous avons vu le nom de Feldbau para\^itre sur la note~\cite{ehresmannfeldbau} le 4 juin; celui de Paul L\'evy para\^it encore le 23 juin sur sa note~\citepyear{LevyNote41}. 
\begin{itemize}
\item Le 4 ao\^ut, puis le 25 ao\^ut, paraissent deux notes du g\'eologue Yaacov Bentor\footnote{Le g\'eologue Yaacov Bentor (1910--2002) sera plus tard un sp\'ecialiste de la Mer Morte et des \'ev\'enements g\'eologiques dont la Bible a gard\'e trace, mais s'int\'eresse ici aux volcans d'Auvergne.}. Alors que la r\`egle est \og Note de M. \textsc{Pr\'enom Nom}, pr\'esent\'ee par M. Pr\'enom Nom\fg et que le pr\'enom en entier est r\'eclam\'e sur de nombreux manuscrits dont les auteurs ont simplement fourni leur initiale\footnote{D'apr\`es Weil~\citepyear[p.~106]{weilsouvenirs}, c'est \`a cette r\`egle que Bourbaki doit d'avoir \'et\'e muni d'un pr\'enom.}, c'est le contraire qui se passe ici, le pr\'enom est ray\'e sur le manuscrit et remplac\'e par son initiale, les deux notes sont publi\'ees sous le nom de Y.~\textsc{Bentor}, moins compromettant.
\item Le 17 novembre, le manuscrit de la note \emph{La survie de souris, de lign\'ee et d'\^age diff\'erents, apr\`es une seule irradiation totale par les rayons $X$}\footnote{Je ne r\'esiste pas au plaisir de mentionner un des r\'esultats de ce travail: \og les animaux de sexe m\^ale se sont montr\'es beaucoup moins r\'esistants \`a l'irradiation totale que les animaux de sexe femelle [sic]\fg.}, de \og M\up{me} \textsc{N.~Dobrolvoska\"ia-Zavadska\"ia}, M. \textsc{S.~V\'er\'etennikoff} et M\up{me} \textsc{M.~Rozd\'evitch}\fg, auteurs auxquels, contrairement \`a la r\`egle \'enonc\'ee ci-dessus et sans doute \`a cause de la longueur de leurs noms, on n'a pas demand\'e d'ajouter leurs pr\'enoms, ce manuscrit est surmont\'e de la mention manuscrite:
\begin{quote}
\noindent Les auteurs sont-ils aryens?
\end{quote}
\item Le 1\up{er} d\'ecembre para\^it finalement, nous l'avons vu, la note~\cite{Ehresmann41}, qui attend depuis le 27 octobre, amput\'ee du nom d'un de ses auteurs.
\end{itemize}
Il semble que fin 1941, tout \'etait r\'egl\'e et que plus aucun math\'ematicien (ni aucun scientifique) juif ne pouvait publier \`a l'Acad\'emie des sciences.

\begin{remarque*}
Le paragraphe pr\'ec\'edent contient \emph{tout} ce que j'ai trouv\'e dans les pochettes de s\'eances. En particulier, je n'ai vu aucune indication que la censure ait v\'erifi\'e que la d\'efinition de \og juif\fg rappel\'ee ici dans la note~\ref{notestatut} s'appliquait. 
\end{remarque*}

\medskip
Quelques scientifiques vont alors penser \`a utiliser le \og pli cachet\'e\fg. Nous verrons Paul L\'evy le faire. Ce sera aussi le cas des biochimistes Jeanne L\'evy et Ernest Kahane\footnote{Jeanne L\'evy (1895--1993) a \'et\'e la premi\`ere femme agr\'eg\'ee de m\'edecine \`a Paris (en 1934). Ernest Kahane (1903--1996) est un chimiste. Ils ont \'ecrit ensemble plusieurs articles et un ouvrage sur le lait. Tous deux \'etaient \og de gauche\fg.} (leur note sur la biochimie de la choline et de ses d\'eriv\'es, parue dans la s\'eance du 6 novembre 1944, \'etait contenue dans un pli cachet\'e d\'epos\'e le 10 juillet 1944).

\section{Andr\'e Bloch, Paul L\'evy, Laurent Schwartz et F\'elix Pollaczek}
\subsection{Le cas d'Andr\'e Bloch (1893--1948)}\label{subsecbloch}
Je renvoie \`a l'article~\cite{CFBloch} pour une biographie de ce math\'ematicien peu ordinaire dont il suffit de rappeler ici qu'apr\`es avoir commis un triple assassinat, il a pass\'e l'essentiel de sa vie professionnelle enferm\'e dans un h\^opital psychiatrique. Bien que semblant de ce fait coup\'e du monde, il est rest\'e en contact avec les math\'ematiques, avec certains math\'ematiciens et, apparemment avec l'actualit\'e. D\`es le 30 d\'ecembre 1940, il \'ecrit \`a \'Emile Picard la lettre suivante (rappelons qu'en d\'ecembre 1940 et encore jusqu'\`a la moiti\'e de l'ann\'ee 1941, l'Acad\'emie des sciences publie des notes d'auteurs juifs, Jacques Feldbau, Laurent Schwartz, Paul L\'evy ou d'autres):

\medskip
\begin{quote}
\begin{flushright}
Saint-Maurice 31 d\'ecembre 1940
\end{flushright}

\medskip
\begin{center}
Monsieur le secr\'etaire perp\'etuel
\end{center}

\medskip
J'ai l'honneur de vous adresser ci-joint une note, dans l'espoir qu'elle vous para\^itra digne d'\^etre accept\'ee \`a l'Acad\'emie.

Intentionnellement, j'ai laiss\'e un blanc apr\`es le titre. Mais, si cela ne pr\'esente pas d'inconv\'enient, vous pouvez le remplacer par mon nom; vous pouvez aussi le remplacer par un nom fictif (Ren\'e Binaud, par exemple) ou m\^eme r\'eel, \`a votre guise.

Je vous serais infiniment reconnaissant, que l'impression de cette note doive \^etre rapide ou ajourn\'ee, de bien vouloir --- tout \`a fait exceptionnellement --- me le faire conna\^itre, par un mot tr\`es court; et je me conformerai \`a l'avenir, en ce qui concerne le sujet de ladite Note, aux indications qui pourront r\'esulter, \`a votre jugement ou ostensiblement, des circonstances dans lesquelles elle aura \'et\'e publi\'ee.

Dans le cas o\`u elle ne pourrait para\^itre, il y aurait lieu de la conserver, sous pli ouvert, dans les archives de l'Acad\'emie.

De toute mani\`ere, je serais tr\`es heureux d'\^etre fix\'e aussi t\^ot que possible, afin de savoir comment je dois organiser ma production scientifique.

Avec mes remerciements anticip\'es, je vous prie d'agr\'eer, Monsieur le Secr\'etaire Perp\'etuel, l'expression de mes sentiments respectueusement et profond\'ement d\'evou\'es.

\medskip
\begin{center}
Andr\'e Bloch, 
\end{center}

\begin{flushright}
57 Grand'Rue\\
Saint-Maurice (Seine)
\end{flushright}

\medskip
PS. Je n'ai pas besoin d'ajouter que, depuis une quinzaine d'ann\'ees et comme tout math\'ematicien s\'erieux, plus que jamais aujourd'hui, je n'\'enonce rien dont je ne sois \`a peu pr\`es s\^ur. Une erreur isol\'ee peut encore m'\'echapper, mais exceptionnellement.
\end{quote}

\medskip
C'est une lettre d\'elicate, allusive, tout \`a fait dans le ton de l'Acad\'emie des sciences, aucun mot excessif n'est explicit\'e, m\^eme l'euph\'emisme \og isra\'elite\fg n'est pas utilis\'e. \`A cette lettre, les secr\'etaires perp\'etuels r\'epondent, le 10 janvier, par une lettre courte, l\'eg\`erement mensong\`ere (Bloch a-t-il vraiment exprim\'e un d\'esir? ce d\'esir-l\`a?), et surtout prudente:

\medskip
\begin{quote}
\begin{center}
Monsieur,
\end{center}

\medskip
Nous avons l'honneur de vous accuser r\'{e}ception de la note que vous venez de nous adresser et qui sera ins\'er\'ee dans les Comptes rendus de l'Acad\'emie.

Suivant le d\'esir que vous exprimez, elle le sera sous le pseudonyme de Ren\'e Binaud.

Veuillez agr\'eer, Monsieur, l'assurance de nos meilleurs sentiments.
\end{quote}

\medskip
La note porte sur une g\'en\'eralisation d'un th\'{e}or\`{e}me de Guldin. Et en effet, Bloch a laiss\'e le nom en blanc sur le manuscrit; une premi\`ere main a rempli le blanc par son vrai nom; une autre a mis \og Ren\'e Binaud\fg; une troisi\`eme main, qui est celle de Picard, a ajout\'e \og pour les CR\fg.

Quelques semaines plus tard, le 20 f\'evrier 1941, nouvelle lettre de Bloch \`a Picard\footnote{La p\'eriode dont il est question ici n'est pas la plus grande p\'eriode cr\'eative de Bloch. Il s'est fait conna\^itre notamment par des articles autour du th\'eor\`eme de Picard et de la r\'epartition des valeurs d'une fonction analytique et semble avoir montr\'e beaucoup de d\'ef\'erence envers \'Emile Picard, qui \'etait \`a l'\'epoque, il faut le dire, un grand pontife de la Science fran\c caise. Il faut aussi noter ici que Picard a toujours \'et\'e un catholique de droite, nationaliste et pratiquant un antis\'emitisme \og ordinaire\fg, comme le montre abondamment sa correspondance avec Lacroix (dossier Picard, archives de l'Acad\'emie des sciences). Citons ici un passage d'une de ses lettres (datant du 29 ao\^ut 1925) o\`u il est question de Bloch: \og Je vois que tout est calme \`a l'Acad\'emie sauf les incartades de notre confr\`ere M.~H. qui est un agit\'e, presque aussi dangereux que son coreligionnaire M.~B. de Charenton\fg (M.~H. d\'esigne Hadamard, M.~B. de Charenton, est Bloch, intern\'e \`a \og Charenton\fg). Un antis\'emitisme ordinaire qui ne l'a pas emp\^ech\'e de publier, en 1941, deux notes de Bloch: c'est un \og agit\'e\fg, un fou dangereux, c'est un juif, mais c'est un math\'ematicien!\label{notepicard}} accompagnant l'envoi d'une nouvelle note. La lettre commence par des remerciements pour la publication de la note Binaud, elle est assez longue et contient beaucoup de politesses, ainsi qu'une discussion de la possible originalit\'e des r\'esultats obtenus (Bloch ne travaille pas, rappelons-le, dans un institut dot\'e d'une biblioth\`eque) qu'il n'est pas indispensable de reproduire ici.

\medskip
\begin{quote}
\begin{flushright}
Saint-Maurice, 20 f\'evrier 1941
\end{flushright}

\medskip
\begin{center}
Illustre Ma\^itre,
\end{center}

\medskip
Veuillez accepter mes plus vifs remerciements pour la lettre tr\`es aimable que vous avez bien voulu m'adresser le mois dernier ainsi que la publication de ma note de g\'eom\'etrie, d'une mani\`ere parfaite, en t\^ete de l'ann\'ee actuelle. [...]

Au sujet de la note que je soumets ci-contre \`a votre jugement, se pose la question de la certitude, et aussi celle de la nouveaut\'e. [...]

En t\^ete de la note, la mention \og Arithm\'etique\fg me para\^it pour plusieurs motifs bien pr\'ef\'erable \`a celle de \og Th\'eorie des nombres\fg, et je suppose que, quoique peu usit\'ee\footnote{Weil raconte dans~\citepyear[p.~116]{weilsouvenirs} qu'\`a Strasbourg dans les ann\'ees 1930, le doyen refusa qu'il intitule un cours \og arithm\'etique\fg parce que \og cela sentait son \'ecole primaire\fg. Les math\'ematiciens utilisaient pourtant ce mot dans son sens actuel au moins depuis le d\'ebut du si\`ecle, me pr\'ecise Catherine Goldstein.}, elle pourra \^etre maintenue.

Je souhaite aussi que le nouveau pseudonyme vous paraisse ad\'equat. Il est en effet pr\'ef\'erable, puisque pseudonyme il y a, que les notes d'arithm\'etique ne paraissent pas sous le m\^eme que celles de g\'eom\'etrie; j'aurais d'ailleurs diff\'erentes choses \`a vous dire \`a ce sujet, mais vous m'excuserez de ne pas insister. \og Marcel Segond\fg ne pr\'esentera, je l'esp\`ere, pas d'inconv\'enient; ce nom propre est rare, et ne figure gu\`ere, je crois, qu'associ\'e \`a un autre. Dans le cas contraire, \og Louis Lechanvre\fg, par exemple, un peu moins idoine, pourra encore convenir.

Enfin, tout en souhaitant ne pas vous importuner, pourrai-je vous rappeler que je vous ai adress\'e il y a deux ans un article \og Sur les propri\'et\'es alg\'ebrico-diff\'erentielles des int\'egrales ab\'eliennes\fg suivi de deux addenda, et dont vous avez bien voulu m'\'ecrire \`a cette \'epoque qu'il serait ins\'er\'e dans votre Bulletin. Peut-\^etre sa publication serait-elle malais\'ee dans les circonstances actuelles, m\^eme sous un pseudonyme; c'est ce dont je ne puis pas me rendre compte; si elle devait avoir lieu, peut-\^etre conviendrait-il pour les \'epreuves de recourir \`a un interm\'ediaire; non, au fait, puisque c'est vous m\^eme ou votre r\'edacteur qui vous en occupez g\'en\'eralement. Quoi qu'il en soit, vous \^etes \`a m\^eme de juger de la question. Si comme il para\^it vraisemblable, sa publication n'\'etait pas possible ou pas d\'esirable pour le moment, je vous serais bien reconnaissant de me le faire retourner et je le garderais pour une \'epoque ult\'erieure; je vous adresserais alors peut-\^etre, peu de temps apr\`es, une note le r\'esumant. Encore une fois, en cette circonstance, je m'en remets enti\`erement \`a votre jugement.

Veuillez agr\'eer, illustre Ma\^itre, avec mes remerciements renouvel\'es, l'assurance de mes sentiments les plus sinc\`erement et respectueusement d\'evou\'es.

\begin{center}
A. Bloch
\end{center}
\end{quote}

\medskip
On aura remarqu\'e que, toujours tr\`es raisonnable, Bloch choisit Segond comme second pseudonyme. La note~\cite{BlochSegond41a} est relue par \'Elie Cartan, accept\'ee, pr\'esent\'ee le 17 mars et publi\'ee dans la s\'eance du 24 mars. Elle porte sur des propri\'et\'es des d\'eveloppements d\'ecimaux des nombres rationnels, une des cons\'equences est le joli r\'esultat suivant (d\^u \`a Kronecker~\citepyear{Kronecker57}): \emph{un entier alg\'ebrique de module inf\'erieur ou \'egal \`a $1$ dont tous les conjugu\'es ont la m\^eme propri\'et\'e est une racine de l'unit\'e}.

L'article dont Bloch attend la parution dans le \emph{Bulletin des Sciences math\'ematiques}, que publiait Picard n'est jamais paru (m\^eme sous pseudonyme). Il est signal\'e comme ayant \'et\'e restitu\'e \`a l'auteur dans la bibliographie de l'article~\cite{CFBloch}. Bloch publie deux autres articles sous le nom de Marcel Segond, tous les deux dans le \emph{Journal de Math\'ematiques pures et appliqu\'ees}, l'un~\citepyear{BlochSegond41} consacr\'e \`a la prolongeabilit\'e des solutions de certaines \'equations diff\'erentielles, le second~\citepyear{BlochSegond42} proche de la note de \og Binaud\fg. Il ne faut pas prendre tr\`es au s\'erieux le d\'edoublement de Binaud-g\'eom\`etre/Segond-arithm\'eticien: il y a des th\'eor\`emes \`a la Guldin sur le volume engendr\'e par une courbe gauche tournant autour d'une droite, dans la note~\cite{Bloch40} (sign\'ee Andr\'e Bloch), dans la note~\cite{BlochBinaud41} (sign\'ee Ren\'e Binaud) et dans l'article~\cite{BlochSegond42} (sign\'e Marcel Segond), ce dernier, \'ecrit en f\'evrier 1940, dit son auteur, faisant m\^eme r\'ef\'erence \`a \og la Note aux Comptes rendus du 27 mai\fg (sans nom d'auteur), c'est-\`a-dire \`a~\cite{Bloch40}, une note sign\'ee du nom d'Andr\'e Bloch. On ne s'attendait sans doute pas \`a ce que les censeurs lisent les articles.

\medskip
La conclusion de cette histoire est conforme \`a son commencement: l'enfermement dans un h\^opital psychiatrique, qui a \'et\'e fatal \`a de nombreux malades et/ou intern\'es pendant l'Occupation nazie (voir par exemple~\cite{Bueltzing}), semble avoir pr\'eserv\'e Bloch des rafles et des pers\'ecutions antis\'emites. Andr\'e Bloch a surv\'ecu \`a l'Occupation\footnote{Il mourra de maladie en 1948. Voir~\cite{CFBloch}.} et, d\'ecid\'ement tr\`es sens\'e, il a publi\'e, d\`es la Lib\'eration (la note para\^it dans la s\'eance du 25 septembre 1944!) sous son vrai nom une note aux \emph{Comptes rendus}~\cite{Bloch44}, dont il n'est pas impossible qu'elle soit le r\'esum\'e de l'article envoy\'e \`a Picard qu'il envisage d'\'ecrire dans la lettre que nous venons de lire. 

\subsection{Le cas de Paul L\'evy (1886--1971)}
La \og strat\'egie\fg de publication (qu'il n'a peut-\^etre pas pens\'ee en ces termes) de Paul L\'evy pendant la p\'eriode concern\'ee montre \`a la fois une vision tr\`es r\'ealiste de la situation et une grande inconscience. Par exemple, juste avant la guerre, \`a une \'epoque o\`u peu de math\'ematiciens fran\c cais publient aux \'Etats-Unis, il envoie un gros article~\citepyear{Levy40} \`a l'\emph{American Journal of Mathematics},
\begin{quote}
\noindent inquiet de l'avenir de l'Europe, j'avais envoy\'e le second aux \'Etats-Unis~\cite[p.~121]{Levy70}.
\end{quote}
Les articles publi\'es pendant la p\'eriode de l'Occupation le sont
\begin{itemize}
\item \cite{Levy40}, \`a l'\emph{American Journal of Mathematics} en 1940, article re\c cu par le journal le 17 octobre 1939,
\item aux \emph{Comptes rendus} en 1941 (puis en 1944), j'ai d\'ej\`a mentionn\'e la note~\cite{LevyNote41} du 23 juin 1941, je reviendrai plus bas sur celle de 1944,
\item \cite{Levy41}, au \emph{Bulletin des Sciences math\'ematiques} en 1941, il s'agit du tout premier fascicule de 1941, qui ne porte pas de date (mais le pr\'ec\'edent vaut pour ao\^ut-septembre-octobre 1940 et le suivant pour mars-avril 1941),
\item \cite{LevySMF41}, au \emph{Bulletin de la Soci\'et\'e math\'ematique de France} en 1941, il s'agit de la publication tardive d'une conf\'erence faite lors d'une r\'eunion de la soci\'et\'e le 6 d\'ecembre 1939,
\item \cite{LevySMFbis41}, aux \emph{Annales de l'Universit\'e de Lyon} en 1941, il s'agit encore d'une conf\'erence faite \`a la \textsc{smf}, section du sud-est, \`a Lyon le 25 janvier 1941 (Paul L\'evy \'etait sans doute \`a Lyon pour donner son cours \`a l'\'Ecole polytechnique qui y avait \'et\'e transf\'er\'ee), qui a \'echapp\'e \`a la publication dans le \emph{Bulletin de la Soci\'et\'e math\'ematique de France} pour une raison que je ne connais pas,
\item \cite{Levy42a, Levy42b}, \`a l'\emph{Enseignement math\'ematique} (deux articles) en 1939--40 (le volume est paru en 1942) dont le deuxi\`eme fait r\'ef\'erence \`a une conf\'erence qui vient d'avoir lieu, en juillet 1939 \`a Gen\`eve, puis dans le volume suivant de ce journal (dat\'e 1942--1950) \cite{Levy42-50},
\item aux \emph{Commentarii Mathematici Helvetici}, un article~\cite{Levy44a} re\c cu le 1\up{er} septembre 1942, puis un autre~\cite{Levy44b} un an plus tard.
\item \`A cette liste d\'ej\`a longue, il convient d'ajouter un article~\cite{Levy43} que Paul L\'evy a envoy\'e le 15 avril 1943 aux \emph{Annales hydrographiques}, publi\'ees par les services hydrographiques de la marine, un endroit inattendu et inhabituel pour un article de math\'ematiques\footnote{Le journal est assez exotique pour que cet article ne soit recens\'e ni par \emph{Math. Reviews}, ni par \emph{Zentralblatt}. J'ai trouv\'e la r\'ef\'erence par hasard dans une liste de ses publications \'etablie par Paul L\'evy et conserv\'ee avec certains de ses papiers dans un dossier par la biblioth\`eque de Chevaleret.\label{hydro} Si le journal n'est pas connu des revues de r\'ef\'erences math\'ematiques, il n'en existe pas moins, et l'article de Paul L\'evy y figure bien.}. Le volume n'est paru qu'en 1946.
\end{itemize}

\medskip
Paul L\'evy a v\'ecu la p\'eriode de l'Occupation dans une semi-clandestinit\'e. Il semble s'\^etre fait quelques illusions, s'\'etonnant par exemple d'apprendre, en novembre 1943, qu'il n'est plus r\'einvesti dans ses fonctions de professeur \`a l'\'Ecole polytechnique depuis le mois d'avril et en particulier qu'on ne lui verse plus de salaire (lettre \`a Fr\'echet\footnote{Maurice Fr\'echet (1878--1973), qui a introduit les espaces m\'etriques en analyse fonctionnelle, est un ami et correspondant de longue date de Paul L\'evy.} du 29 novembre 1943~\cite{LevyFrechet}). Il poss\'edait des faux papiers au nom de Paul Leng\'e, se faisait \'ecrire sous le nom de Paul Piron chez son gendre Robert Piron \`a Grenoble, mais n'a pas essay\'e d'utiliser un pseudonyme pour publier ses articles.

Il a r\'eussi \`a garder des contacts \`a Gen\`eve, qui lui ont permis de publier deux articles \`a \emph{l'Enseignement math\'ematique} en 1942 et deux autres \`a \emph{Commentarii} en 1943 et 1944.

Le 1\up{er} juin 1943, dans une lettre \`a Fr\'echet~\cite{LevyFrechet}, il \'ecrit (confirmant l'interdiction de publication aux \emph{Comptes rendus}):

\medskip
\begin{quote}
Le r\'esultat obtenu me para\^it assez important et, comme je ne peux pas pr\'esenter en ce moment de Note \`a l'Acad\'emie, je vous serais oblig\'e de conserver la pr\'esente lettre.

Il s'agit de la fonction al\'eatoire du mouvement brownien, $X(t)$, d\'ej\`a \'etudi\'ee dans mon livre sur les variables al\'eatoires et dans mes deux m\'emoires de 1939, ainsi que par divers savants \'etrangers.

[...]

Finalement [...] on voit que

\underline{Th\'eor\`eme}: \emph{La nature stochastique de $E$ est invariante par n'importe quelle homographie.}

[...]
\end{quote}

Il s'agit d'un r\'esultat que Paul L\'evy va utiliser dans un de ses articles suisses~\citepyear{Levy44b} et dont il a l'id\'ee brillante de le r\'ediger et de l'envoyer \`a l'Acad\'emie... mais comme un pli cachet\'e, qu'aucune loi n'emp\^eche l'Acad\'emie de recevoir le 16 juin 1943, d'enregistrer sous le num\'ero 11904 et de conserver jusqu'au retour de jours meilleurs. Et en effet, le pli est ouvert d\`es le 23 octobre 1944 et, imm\'ediatement publi\'e, il devient la note~\citepyear{LevyNote44}.

Quant \`a l'article suisse~\citepyear{Levy44b} dans lequel le r\'esultat est utilis\'e, Paul L\'evy l'a envoy\'e \`a Wavre pour les \emph{Commentarii} le 23 ao\^ut 1943, apr\`es qu'il ait \'et\'e refus\'e par les \emph{Annales de l'Universit\'e de Grenoble}. Il peut sembler \'etrange, et surtout \'etonnamment imprudent d'avoir envoy\'e un article \`a ce journal en 1943, mais nous verrons que son gendre Laurent Schwartz en faisait autant \`a Toulouse --- qui \'etait peut-\^etre moins expos\'ee. 

\begin{remarque*}[Sur les Annales de l'Universit\'e de Grenoble]
C'est un journal g\'en\'era\-liste (sciences et m\'edecine). Ses volumes 19 (1943) et 20 (1944) --- ceux qui auraient pu accueillir~\cite{Levy44b} --- ne contiennent chacun qu'un article de math\'ematiques --- dans les deux cas, il s'agit d'un article de Brelot\footnote{Marcel Brelot (1903--1987) \'etait professeur \`a Grenoble et sera le v\'eritable cr\'eateur des \emph{Annales de l'Institut Fourier} (voir~\cite{ChoquetBrelot}).}.

La r\'egion de Grenoble et avec elle son universit\'e, semblent en effet avoir \'et\'e serr\'ees de pr\`es par les forces d'occupation, \`a cause des maquis de l'Is\`ere et des activit\'es de la r\'esistance dans la r\'egion, la ville et l'universit\'e\footnote{La r\'ealit\'e de la r\'esistance et de la collaboration dans cette institution est \'evidemment nuanc\'ee, voir~\cite{Dereymez}.} ainsi que du grand nombre de r\'efugi\'es que les montagnes accueillaient. Dans la lettre du 29 novembre 1943 d\'ej\`a cit\'ee~\cite[p.~212]{LevyFrechet}, Paul L\'evy explique \`a Fr\'echet:

\medskip
\begin{quote}
Apr\`es trois mois de vie errante, je suis \`a nouveau install\'e dans une maison o\`u j'esp\`ere passer l'hiver\footnote{Il n'y a pas sur cette lettre d'indication du lieu d'o\`u elle a \'et\'e \'ecrite. Pour une fois, son auteur s'est montr\'e prudent.}. Grenoble, comme vous devez le savoir, est tr\`es agit\'e; il y a constamment des incidents: attentats, sanctions, repr\'esailles. Ma famille a insist\'e il y a trois mois pour que je ne reste pas au voisinage d'une ville si expos\'ee, d'autant plus qu'\'etant connu comme je l'\'etais, je pouvais \^etre particuli\`erement expos\'e.
\end{quote}

\medskip
La ville de Grenoble recevra la Croix de la Lib\'eration, des mains de de Gaulle, le 5 novembre 1944. Il est assez \'emouvant, apr\`es avoir lu les lettres de Paul L\'evy, d'ouvrir le tout premier volume des \emph{Annales de l'Universit\'e de Grenoble} qui para\^it apr\`es la Lib\'eration. Pas seulement parce qu'il commence par une s\'erie d'articles n\'ecrologiques sur les cadres de l'universit\'e morts pendant l'Occupation, mais surtout parce qu'il continue par des articles scientifiques, dont le premier est d\^u \`a Paul L\'evy~\citepyear{Levy45} et a \'et\'e r\'edig\'e, en 1944, \emph{\`a la demande} du Directeur de l'Institut polytechnique de Grenoble (c'est-\`a-dire F\'elix Esclangon, d\'ej\`a \`a ce poste depuis 1941 --- comme les autres journaux, celui-l\`a ne contient aucune mention des noms de ses r\'edacteurs)... et dont le deuxi\`eme est d\^u \`a Laurent Schwartz~\citepyear{Schwartz45}.
\end{remarque*}

\begin{remarque*}[Sur le Bulletin de la Soci\'et\'e math\'ematique de France]
\`A la fois journal scientifique et organe de la soci\'et\'e, le \emph{Bulletin} publie des articles de recherche et des informations sur la vie de la soci\'et\'e (liste des membres avec leurs adresses, \'election de nouveaux membres, composition du Conseil et du Bureau, et m\^eme parfois des expos\'es de math\'ematiques pr\'esent\'es \`a telle ou telle de ses r\'eunions).

Pendant la p\'eriode concern\'ee, la liste des membres de la soci\'et\'e dispara\^it (prudence \'el\'ementaire). Pas celle des membres du conseil. Quant \`a la r\'edaction, ce sont, tout simplement, les secr\'etaires qui r\'edigent. Bref, le \emph{Bulletin de la Soci\'et\'e math\'ematique de France}, pendant l'Occupation, c'est Henri Cartan.

En 1941, le \emph{Bulletin} publie avec un peu de retard (mais sans doute juste \`a temps) des conf\'erences donn\'ees par, notamment, Paul L\'evy et Szolem Mandelbrojt\footnote{Szolem Mandelbrojt (1899--1983) a pass\'e la dur\'ee de la guerre aux \'Etats-Unis, nous l'avons dit.}. 

Retour \`a Feldbau-Laboureur. Henri Cartan est un ami de Jacques Feldbau, qu'il a connu \`a Strasbourg avant la guerre. Il va l'aider, personnellement, lui et aussi sa famille, pendant toute la dur\'ee de la guerre. L'affaire de la disparition de Feldbau de~\cite{Ehresmann41} a certainement \'et\'e discut\'ee par \'Elie Cartan avec son fils. L'id\'ee de publier sous un pseudonyme, comme \'Elie Cartan savait tr\`es bien que Bloch le faisait, sans doute aussi (d'autant plus qu'Henri Cartan avait commenc\'e ses travaux math\'ematiques autour d'une conjecture d'Andr\'e Bloch). M\^eme si ce n'est pas lui qui l'a sugg\'er\'e, Henri Cartan savait tr\`es bien qui \'etait Jacques \og Laboureur\fg lorsqu'il a publi\'e les deux notes~\cite{Laboureur42, Laboureur43} dans \og son\fg journal. C'est d'ailleurs de ce nom que Feldbau signera le mot qu'il lui enverra pour lui demander de prendre soin de sa famille lorsqu'il quittera Drancy pour Auschwitz\footnote{Henri Cartan a r\'eussi \`a envoyer des colis \`a Feldbau \`a Auschwitz. Lorsque les survivants des camps ont commenc\'e \`a revenir en 1945, il a fait tout son possible pour collecter des informations et savoir ce qu'il \'etait advenu de Feldbau. Voir~\cite{feldbau} pour les sources de ces informations.} (voir~\cite{georgescerf}) en octobre 1943.
\end{remarque*}

\subsection{Le cas de Laurent Schwartz (1915--2002)}
Comme Jacques Feldbau, Laurent Schwartz est alors un jeune math\'ematicien, il est n\'e en 1915, a fait son service militaire avant la guerre, puis a fait la guerre. Enfin il se met \`a la recherche pendant l'Occupation. Il publie alors
\begin{itemize}
\item deux notes~\citepyear{Schwartz41a,Schwartz41c} dans le volume du premier semestre 1941 des \emph{Comptes rendus},
\item sa th\`ese~\citepyear{Schwartz43a}, soutenue le 9 janvier 1943 \`a Clermont-Ferrand et publi\'ee par les \'editions Hermann, dont le directeur, Freymann, n'h\'esitait pas, non seulement \`a publier la th\`ese d'un math\'ematicien juif \`a Paris en 1943, mais m\^eme \`a afficher dans sa vitrine les \oe uvres d'Einstein (voir ce que raconte Schwartz lui-m\^eme dans son livre~\citepyear[p.~175]{schwartz}),
\item un chapitre suppl\'ementaire~\citepyear{Schwartz43b} qu'il n'a pas inclus dans celle-ci et que publient les \emph{Annales de l'Universit\'e de Toulouse}, en 1943, elles aussi. Dans son livre~\cite[p.~174]{schwartz}, Schwartz signale cet \og article compl\'ementaire\fg sans autre commentaire.
\end{itemize}

\`A la Lib\'eration, il ne perdra pas beaucoup de temps et enverra rapidement un court article~\citepyear{Schwartz44} au \emph{Bulletin de la Soci\'et\'e math\'ematique de France}, que celui-ci recevra le 31 octobre 1944. J'ai mentionn\'e ci-dessus un article ult\'erieur~\citepyear{Schwartz45}.

\begin{remarque*}[Sur les Annales de l'Universit\'e de Toulouse]
Il semble qu'ici, comme dans le cas du \emph{Bulletin} et d'Henri Cartan, le r\^ole du secr\'etaire a d\^u \^etre d\'eterminant. Le secr\'etaire des \emph{Annales de Toulouse}, c'est Adolphe Buhl\footnote{Adolphe Buhl (1878--1949) \'etait professeur \`a l'universit\'e de Toulouse et dirigeait aussi, avec H. Fehr, l'\emph{Enseignement math\'ematique} (voir~\cite{FehrBuhl}).}. Tranquillement, Adolphe Buhl publie l'article~\citepyear{Schwartz43b} de Schwartz dans son journal. Celui-ci \'evoque les difficult\'es de l'heure et le remercie:

\begin{quote}
Comme je l'ai d\'ej\`a dit dans l'introduction de ma th\`ese, et pour les m\^emes raisons, je m'excuse \`a l'avance des lacunes ou des erreurs bibliographiques, cons\'equences d'une documentation rendue difficile par les circonstances.

Je tiens \`a exprimer toute ma reconnaissance \`a M. Buhl, qui s'est occup\'e de l'insertion de mon M\'emoire dans les pr\'esentes Annales.
\end{quote}

Adolphe Buhl \'etait r\'eput\'e un homme courageux et droit (voir~\cite{FehrBuhl}). Apr\`es tout, nous l'avons dit, la loi fran\c caise n'interdisait pas aux juifs de publier des articles scientifiques. Il est probable aussi que Toulouse \'etait moins expos\'ee que Grenoble. Peut-\^etre aussi Schwartz \'etait-il moins visible que Paul L\'evy. Si le nom de Paul L\'evy avait sans mal \'et\'e identifi\'e comme juif, ce n'\'etait peut-\^etre pas le cas de celui de son gendre Laurent Schwartz, qui n'\'etait pas aussi connu \`a l'\'epoque qu'il ne l'est aujourd'hui.
\end{remarque*}

\subsection{Le cas de F\'elix Pollaczek (1892--1981)}\label{subsecpollaczek}
F\'elix Pollaczek, n\'e \`a Vienne en 1892, avait perdu son emploi \`a Berlin et quitt\'e l'Allemagne apr\`es l'arriv\'ee des nazis au pouvoir. Il a ensuite v\'ecu en France, o\`u il a notamment pass\'e la p\'eriode de l'Occupation, mais je ne sais pr\'ecis\'ement ni o\`u ni dans quelles conditions (voir le court article n\'ecrologique~\cite{CohenPollaczek}). Toujours est-il qu'il a publi\'e un article de math\'ematiques~\citepyear{Pollaczek42} dans les \emph{Annales de l'Universit\'e de Lyon} en 1942 (F\'elix Pollaczek est un probabiliste bien connu) et trois articles~\citepyear{Pollaczek43a,Pollaczek43b,Pollaczek43c} de physique (il \'etait, avant la guerre, ing\'enieur consultant pour la \emph{Soci\'et\'e d'\'etudes pour les liaisons t\'el\'ephoniques et t\'el\'egraphique}) dans les \emph{Cahiers de physique} en 1943.

\begin{remarque*}[Sur les Annales de l'Universit\'e de Lyon]
Dans les \emph{Annales de l'Universit\'e de Lyon}, cit\'ees ci-dessus, on a reconnu le journal ayant accueilli l'article~\cite{LevySMFbis41}. Je ne sais pas qui \'etaient les r\'edacteurs de ce journal, dont la section A publiait des articles de math\'ematiques et d'astronomie, mais l'\'editeur commercial \'etait Hermann. Outre~\cite{LevySMFbis41} et~\cite{Pollaczek42}, notons la pr\'esence dans ce journal du premier article~\cite{Samuel42} du tout jeune Pierre Samuel (n\'e en 1921) dont l'entr\'ee \`a l'\textsc{ens} avait \'et\'e diff\'er\'ee apr\`es sa r\'eussite au concours... en 1940.
\end{remarque*}

\begin{remarque*}[Sur les Cahiers de physique]
Les \emph{Cahiers de physique} avaient \'et\'e fond\'es \`a Marseille en 1941 par Charles Fabry (un des physiciens acad\'emiciens r\'eput\'e \og de gauche\fg d\'ej\`a mentionn\'e ci-dessus et un ami d'Henri Abraham), pour permettre aux physiciens de la zone libre de publier leurs travaux. Parmi les noms des auteurs dont ce journal a accept\'e et fait para\^itre des articles pendant l'Occupation, on aussi peut relever celui du physicien Th\'eo Kahan, qui aurait sans doute fait tiquer le secr\'etaire de l'Acad\'emie des sciences. Ce qui semble confirmer qu'il n'\'etait pas impossible (au moins en zone sud, mais m\^eme apr\`es novembre 1942) \`a un r\'edacteur de publier les articles qu'il avait d\'ecid\'e de publier.
\end{remarque*}

\section*{Conclusions et perspectives}
Les (peu nombreux) math\'ematiciens fran\c cais qui ont v\'ecu la p\'eriode de l'Occupation en France et \`a qui une application tr\`es rigoureuse des lois antis\'emites de l'\og\'Etat fran\c cais\fg interdisait de faire para\^itre leurs r\'esultats (\`a l'Acad\'emie des sciences notamment), ont finalement trouv\'e des journaux qui ont accept\'e leurs articles. En affinant cette remarque, on note que les journaux de la zone sud dont il a \'et\'e question (\emph{Annales de l'Universit\'e de Toulouse}, \emph{Annales de l'Universit\'e de Lyon}, \emph{Cahiers de physique}\footnote{Nous avons not\'e aussi la d\'etermination d'un groupe de physiciens \`a ne pas s'accommoder des lois antis\'emites.} \`a Marseille) ont publi\'e les articles de ces auteurs sous leur v\'eritable identit\'e, alors que le \emph{Bulletin de la \textsc{smf}}, publi\'e \`a Paris, lorsqu'il l'a fait, l'a fait sous pseudonyme.

Ces mani\`eres anormales de publier a eu pour effet d'emp\^echer le jeune topologue Jacques Feldbau, mort en d\'eportation et qui donc n'\'etait plus l\`a pour continuer \`a faire avancer la th\'eorie, d'\^etre reconnu pour son apport \`a la topologie comme il aurait d\^u l'\^etre. Il s'agissait d'un jeune math\'ematicien essayant de publier ses premiers r\'esultats dans un domaine nouveau, en pleine cr\'eation, avec une concurrence s\'erieuse de coll\`egues suisses, anglais et am\'ericains qui travaillaient, publiaient (et, tout simplement, vivaient) dans des conditions plus faciles.

\medskip
D'autres math\'ematiciens fran\c cais n'ont pas \'eprouv\'e les m\^emes difficult\'es \`a faire conna\^itre leurs travaux (ni m\^eme \`a vivre). Certains ont v\'ecu cette p\'eriode aux \'Etats-Unis (Chevalley, Hadamard, Mandelbrojt, Weil) et ont naturellement publi\'e dans des journaux am\'ericains\footnote{Dans \emph{Amer. J. of math.} (Chevalley), \emph{Ann. of math.} (Chevalley, Hadamard), \emph{Bull. Amer. Math. Soc.} (Chevalley, Hadamard), \emph{Duke Math. J.} (Mandelbrojt), \emph{Trans. Amer. Math. Soc.} (Chevalley, Mandelbrojt, Weil).}. D'autres, bien que vivant en France, ont pu continuer \`a publier plus ou moins normalement, malgr\'e les difficult\'es de la vie quotidienne. Parmi eux, plusieurs ont publi\'e dans des journaux allemands, certains ont trouv\'e des \og accommodements\fg avec la situation politique, quelques-uns choisissant m\^eme la collaboration. Le cas des prisonniers de guerre ou, plus exactement, les cas vari\'es des prisonniers de guerre (Leray, Pauc, Roger, Ville, ou d'autres) est (sont) aussi digne(s) d'int\'er\^et.

En pr\'eparant cet article, j'ai bien entendu trouv\'e des mat\'eriaux sur ces questions aussi. Il est donc vraisemblable que j'en \'etudierai certaines ult\'erieurement. Les difficult\'es mat\'erielles et th\'eoriques sont beaucoup plus grandes et ceci pour au moins deux raisons pas compl\`etement disjointes:
\begin{itemize}
\item la chape de silence\footnote{Si les deux livres de souvenirs~\cite{schwartz,weilsouvenirs} de math\'ematiciens t\'emoins de cette \'epoque dont nous disposons nous ont donn\'e peu de renseignements sur les questions envisag\'ees dans le pr\'esent article, ils en donnent encore moins sur ces autres questions. Il n'y a par exemple que deux allusions \`a un (m\^eme) math\'ematicien collaborateur dans~\cite[p.~152]{schwartz}, mais il n'est pas nomm\'e (m\^eme si identifiable).} que la communaut\'e math\'ematique elle-m\^eme a d\'epos\'ee sur ces questions d\`es apr\`es la Lib\'eration en a fait dispara\^itre une bonne partie de la m\'emoire (en m\^eme temps que de nombreux documents ont sans doute \'et\'e d\'etruits); elle est peut-\^etre encore trop lourde pour \^etre soulev\'ee,
\item car en effet, il n'est pas s\^ur qu'il soit encore aujourd'hui tr\`es facile de parler, en historien$\cdot$ne, de fa\c con rigoureuse et sereine, de la Collaboration --- \`a propos de coll\`egues, de membres de notre petit monde, et pas de criminels de guerre!
\end{itemize}

\medskip
Je conclurai sur cette derni\`ere question, en revenant \`a la note~\ref{notepicard} pour l'illustrer par l'exemple du math\'ematicien \'Emile Picard, catholique r\'eactionnaire, profond\'ement anti-allemand, p\'etainiste\footnote{Lettre \`a Lacroix, 11 juillet 1940, dossier Picard, archives de l'Acad\'emie des sciences: \og J'admire le Mar\'echal P\'etain\fg.}, presque collaborationniste\footnote{Lettre \`a Lacroix, 30 janvier 1941 (exactement contemporaine de la publication de~\cite{BlochBinaud41}), dossier Picard, archives de l'Acad\'emie des sciences: \og Pour le fond, je suis d'accord avec Claude qu'une collaboration tr\`es g\'en\'erale avec l'Allemagne est n\'ecessaire sous peine d'un \'{e}crasement complet de la France pour un temps ind\'efini; collaboration accept\'ee en principe par le Mar\'echal \`a Montoire [...] Mais ce qu'on peut reprocher \`a Claude c'est la confusion entre la collaboration \emph{scientifique} [...] et la collaboration \'economique et politique [...]\fg (il est question ici de Georges Claude, unique membre de l'Acad\'emie des sciences dont l'\'election sera annul\'ee \`a la Lib\'eration pour faits de collaboration).}, dont un des r\'esultats concrets de l'accommodement avec les lois antis\'emites est la publication des notes d'Andr\'e Bloch, au tout d\'ebut de 1941 il est vrai, mais Picard est mort le 11 d\'ecembre 1941 (au moment de la disparition de Feldbau de la note~\cite{Ehresmann41})...

\subsection*{Remerciements}

Pour ce travail, j'ai d\^u consulter de nombreux articles publi\'es dans les ann\'ees 1930 et 40, ainsi que d'autres parus plus r\'ecemment dans les \emph{Cahiers du S\'eminaire d'histoire des math\'ematiques}. Je tiens \`a rendre hommage ici \`a l'exceptionnelle qualit\'e du travail de num\'erisation r\'ealis\'e par l'\textsc{ums} \textsc{MathDoc} dans son programme \textsc{numdam}, qui a rendu cette consultation facile et agr\'eable.

\medskip
Je remercie le service des archives de l'Acad\'emie des sciences, en les personnes de Florence Greffe, Claudine Pouret et surtout Christiane Pavel et Pierre Leroi pour la gentillesse avec laquelle ils ont cherch\'e et sorti, puis rang\'e la bonne centaine de pochettes que j'ai ouvertes pour r\'ealiser ce travail. Je remercie aussi l'Acad\'emie des sciences pour l'autorisation de publier des extraits des deux lettres d'Andr\'e Bloch et de celle d'\'Emile Borel\footnote{N'ayant pas connaissance d'\'eventuels ayants droit, ni d'Andr\'e Bloch ni d'\'Emile Borel, je n'ai pas demand\'e d'autre autorisation.}.

\medskip
Toute ma reconnaissance aussi \`a, par ordre alphab\'etique,

\begin{itemize}
\item Catherine Goldstein pour ses critiques amicales d'une premi\`ere version de ce texte, ses suggestions et aussi pour les r\'ef\'erences qu'elle m'a donn\'ees et les informations qu'elle m'a aid\'ee \`a trouver,

\item Michel Pinault\footnote{Michel Pinault est l'auteur du livre~\cite{pinaultlivre} sur Fr\'ed\'eric Joliot-Curie, indispensable pour la compr\'ehension du contexte des milieux scientifiques pendant l'\'epoque concern\'ee.}, qui m'a accompagn\'ee et guid\'ee lors de mes premi\`eres visites aux archives de l'Acad\'emie des sciences, pour l'aide qu'il m'a apport\'ee, les informations qu'il m'a donn\'ees et les documents qu'il m'a communiqu\'es, la lettre du \textsc{mbf} \`a Fourneau reproduite au~\T\ref{secacad} notamment,

\item Norbert Schappacher pour son aide, notamment avec les lettres de Hasse, et pour ses remarques sur une version pr\'eliminaire de cet article,

\item Michel Zisman, surtout, dont l'article~\citepyear{zismanjames} est l'une des origines de ce travail et dont les encouragements m'ont \'et\'e pr\'ecieux,

\item Liliane Zweig pour son aide avec le carton d'archives de Paul L\'evy dont il est question dans la note~\ref{hydro}.
\end{itemize}

\medskip
Enfin je remercie deux \og anonymes\fg qui ont \'ecrit pour la \textsc{rhm} des rapports sur une version pr\'ec\'edente de cet article, le math\'ematicien pour la pertinence de ses remarques et pour ses utiles propositions d'am\'elioration et l'historien des math\'ematiques pour sa suggestion de consid\'erer le cas de F\'elix Pollaczek.

\nocite{facssousvichy}

\newcommand{\SortNoop}[1]{}

\vfill

\end{document}